\documentclass[12pt]{article}

\usepackage{float}
\usepackage[margin=1in]{geometry}
\usepackage[T1]{fontenc}
\usepackage[utf8]{inputenc}
\usepackage{lmodern}
\usepackage{setspace}
\usepackage{amsmath,amssymb,amsthm,mathtools}
\usepackage{bm}
\usepackage{graphicx}
\usepackage{xcolor}
\usepackage{booktabs}
\usepackage{array}
\usepackage{enumitem}
\usepackage{caption}
\usepackage{subcaption}
\usepackage{hyperref}
\usepackage[nameinlink,noabbrev]{cleveref}
\usepackage[numbers,sort&compress]{natbib}

\hypersetup{
    colorlinks=true,
    linkcolor=blue,
    citecolor=blue,
    urlcolor=blue,
    pdftitle={Feature-Based Continuation of Pattern Transitions in a One-Dimensional Brusselator},
    pdfauthor={Qiushi Yu}
}

\theoremstyle{definition}

\title{%
    \textbf{Feature-Based Continuation of Pattern Transitions in a One-Dimensional Brusselator}
}
\author{Qiushi Yu}
\date{August 2026}

\begin{document}

\maketitle

\begin{abstract}
Different long-time patterns can prevail in different regions of a reaction--diffusion system's parameter space. We study the transition curves between such regimes for a one-dimensional Brusselator in the two-parameter plane $(\sigma,b)$, focusing on wave/stripe-like, spiral/source-defect-like, and target-like states. We develop a feature-based continuation framework built on time-dependent PDE simulations. Scalar observables extracted from late-time solution data distinguish the regimes and define threshold level sets where the feature crossings are regular. A secant predictor and a local one-dimensional sweep corrector are used to trace these level sets. For the spiral transition, we introduce a branch-adapted spacetime symmetry-defect feature that separates asymmetric source-like patterns from more symmetric wave patterns. For the target transition, we use the minimum of a core spatial-variance score and a tail temporal-variance score to detect the characteristic structure of half-target states. The method recovers robust side portions of both transition curves. In lower parameter regions, where mixed and irregular patterns make a single scalar feature less specific, we instead report transition estimates obtained from vertical parameter sweeps and direct inspection of spacetime plots. These results show how simulation-based continuation and direct pattern classification can be combined to map regime boundaries while preserving the different levels of numerical evidence.
\end{abstract}

\section{Introduction}

\subsection{Pattern Formation in Reaction-Diffusion Systems}
Reaction-diffusion systems form a broad and important class of partial differential equations. In their simplest form, such models describe the evolution of concentrations, densities, or other distributed quantities through the interplay of local production--depletion mechanisms and spatial transport. This interaction can generate dynamics that are not present in either purely kinetic or purely diffusive models alone. As a result, reaction-diffusion equations have become a central framework for studying self-organization in chemical, physical, and biological systems \cite{CrossHohenberg1993,Murray2003,VanagEpstein2009}. Early chemical models of symmetry breaking in dissipative media already illustrated how nonequilibrium kinetics can generate spatial structure, helping to establish the broader conceptual background for later reaction--diffusion pattern theory \cite{PrigogineLefever1968}.

One of the most important features of reaction-diffusion systems is their ability to produce coherent spatial and spatiotemporal patterns \cite{CrossHohenberg1993,EpsteinPojman1998,Hoyle2006}. Depending on the model and parameter regime, one may observe homogeneous steady states, traveling waves or wave trains, stripes, spots, target patterns, spiral waves, and localized structures such as defects or interfaces between competing states \cite{CrossHohenberg1993,Stich2003,ZhaoMaffaSandstede2025}. These patterns arise in a wide range of settings, from chemical reactions and fluid systems to biological media and other nonequilibrium processes \cite{CrossHohenberg1993,Murray2003,ZhaoMaffaSandstede2025}. In particular, reaction-diffusion theory provides a natural language for describing how microscopic interactions give rise to macroscopic pattern selection \cite{CrossHohenberg1993,Hoyle2006}.

An important point is that these structures are typically \emph{far-from-equilibrium} phenomena \cite{CrossHohenberg1993,EpsteinPojman1998}. Unlike equilibrium configurations determined by free-energy minimization, dissipative patterns require ongoing interaction between kinetics, transport, and often external forcing or sustained fluxes \cite{CrossHohenberg1993,VanagEpstein2009}. As a result, the system can support multiple qualitatively distinct long-time behaviors, and these behaviors may change as parameters vary \cite{CrossHohenberg1993,Stich2003,VanagEpstein2009}. Even in relatively simple models, modifying a control parameter can change the dominant observed pattern, for example from a spatially uniform state to stripes or spots \cite{CrossHohenberg1993,Pearson1993,TzouMaBaylissMatkowskyVolpert2013}.

For this reason, it is often not enough to study individual solutions in isolation. From a dynamical-systems viewpoint, it is important to understand how parameter space is partitioned into regions with different prevailing long-time regimes \cite{CrossHohenberg1993,ZhaoMaffaSandstede2025}. To address this problem is to determine the boundaries between parameter regions where different patterns dominate \cite{ZhaoMaffaSandstede2025}. Zhao, Maffa, and Sandstede make this perspective explicit: rather than focusing only on a single coherent structure, one can ask which pattern is prevalent in a given part of parameter space and then attempt to trace the corresponding bifurcation or transition curves separating those regions \cite{ZhaoMaffaSandstede2025}. This viewpoint is especially relevant in systems where direct simulation reveals several different classes of patterns, but where the location of the transition between them is not immediately clear \cite{ZhaoMaffaSandstede2025}.

In the present work, the starting point for the study of a one-dimensional Brusselator model is precisely this general perspective. The main goal is to understand how different regimes are organized in parameter space and how one can compute the transition curves between them, rather than simply exhibit particular solutions \cite{ZhaoMaffaSandstede2025}. Before introducing the specific continuation framework used here, it is therefore essential to place the problem in the broader context of pattern formation in reaction-diffusion systems \cite{CrossHohenberg1993,Murray2003,Hoyle2006}.

\subsection{Defects, sources, and one-dimensional analogues}
In oscillatory media, the relevant long-time behavior is often not described solely by stationary spatial patterns. Instead, one frequently encounters \emph{wave trains}, that is, spatially periodic traveling waves that propagate through the medium while maintaining a coherent profile \cite{CrossHohenberg1993,SandstedeNotes,DoelmanSandstedeScheelSchneider2009}. Such wave trains provide a natural asymptotic background for many more localized structures, and they play a central role in the dynamical-systems description of oscillatory pattern formation \cite{CrossHohenberg1993,Stich2003,SandstedeScheel2004}. In this setting, the important objects are not only domain-filling patterns, but also localized structures that connect, emit, or reorganize wave trains \cite{SandstedeNotes,SandstedeScheel2004}.

A convenient framework for describing such localized objects is the notion of a \emph{defect}. This emphasis on localized organizing structures is also consistent with earlier work on localized patterns in reaction--diffusion systems, where spatially confined structures were already recognized as dynamically significant objects rather than mere transient irregularities \cite{KogaKuramoto1980}. In the dynamical-systems framework developed by Sandstede and Scheel, defects are coherent structures that are asymptotic to wave trains and can be classified according to the relation between defect speed and the group velocities of the asymptotic waves. \cite{SandstedeNotes}. Thus, a defect is not merely an isolated irregularity in a pattern, but rather a dynamically organized interface-like structure embedded in an oscillatory background \cite{SandstedeNotes,SandstedeScheel2004}. This idea is especially useful because it connects localized wave-emitting or wave-selecting behavior with the broader theory of coherent structures in partial differential equations \cite{ChampneysSandstede2007,SandstedeScheel2004}.

Among the different classes of defects, \emph{sources} are especially relevant for the present work. A source defect acts as an organizing center, or pacemaker, for the surrounding medium by emitting waves outward \cite{SandstedeNotes,Stich2003}. In the language of oscillatory media, such objects generate wave trains to both sides and select their asymptotic wavenumbers and propagation behavior \cite{SandstedeNotes,SandstedeScheel2004}. Sandstede and Scheel's classification places sources alongside sinks, contact defects, and transmission defects, with the classification determined by comparing the defect speed with the group velocities of the asymptotic wave trains \cite{SandstedeNotes,SandstedeScheel2004}. This classification is conceptually important because it shows that a localized wave-emitting state is part of a systematic dynamical taxonomy rather than an isolated special case \cite{SandstedeScheel2004}.

This language remains highly relevant even in one spatial dimension. Although spiral waves are usually understood as genuinely two-dimensional objects, one-dimensional systems can still support source-like or defect-like structures that organize outward-propagating waves \cite{PerraudBorckmans1993,CytrynbaumLewis2009}. In such cases, the geometry is simpler, but the essential dynamical role remains the same: a localized core mediates and selects the wave behavior observed in the far field \cite{PerraudBorckmans1993,SandstedeScheel2004,CytrynbaumLewis2009}. For the purposes of the present work, this is precisely the point. The patterns of interest are not classical two-dimensional spirals themselves, but rather one-dimensional source-like or spiral-like states that separate or organize different oscillatory regimes \cite{PerraudBorckmans1993,CytrynbaumLewis2009}.

The one-dimensional setting is therefore neither artificial nor trivial. On the one hand, it captures essential source and defect dynamics in a form that is mathematically and computationally accessible \cite{PerraudBorckmans1993,CytrynbaumLewis2009}. On the other hand, it provides a natural arena for studying regime transitions and for developing numerical continuation methods that track the boundaries between different long-time behaviors \cite{ZhaoMaffaSandstede2025,ChampneysSandstede2007,Uecker2021}. In particular, the reduced geometric complexity of one dimension makes it possible to focus more directly on the detection and selection of wave-emitting structures, while still retaining the core dynamical phenomena that motivate the broader defect framework \cite{SandstedeScheel2004,CytrynbaumLewis2009}.

For these reasons, the one-dimensional Brusselator problem studied in this work is naturally viewed through the lens of defects and sources. This perspective clarifies why localized wave-emitting states matter, why transitions between such states and more symmetric wave regimes are dynamically meaningful, and why a feature-based continuation framework can be interpreted as a tool for tracing boundaries between distinct classes of defect-mediated behavior \cite{SandstedeScheel2004,ZhaoMaffaSandstede2025}.

\subsection{Relation to Target Patterns and One-Dimensional ``Spirals''}

The pattern regimes studied in this work are closely related to two broader strands of literature: the theory of target patterns and pacemakers, and the study of one-dimensional or quasi-one-dimensional asynchronous wave sources \cite{Stich2003,PerraudBorckmans1993}. These two viewpoints provide a natural conceptual bridge between the defect framework discussed above and the specific spiral-like and target-like states observed in the one-dimensional Brusselator.

A \emph{target pattern} consists of outward-propagating waves emitted by a localized wave source, often called a \emph{pacemaker} \cite{KopellHoward1981,Stich2003}. In physical space, the classical two-dimensional picture is that of concentric waves radiating from a core, but the essential dynamical feature is more general: a localized structure entrains the surrounding medium by repeatedly generating outgoing waves \cite{Stich2003,SandstedeScheel2021}. Stich emphasizes this broader viewpoint, describing target patterns and pacemakers as generic wave-source structures in reaction--diffusion systems and distinguishing between heterogeneous pacemakers, which arise from imposed inhomogeneities, and self-organized pacemakers, which emerge intrinsically from the dynamics of the medium itself \cite{Stich2003}. This pacemaker literature is important here because it broadens the interpretation of wave-emitting structures beyond the geometry of classical spiral waves alone \cite{Stich2003,KopellHoward1981}.

The relevance of this viewpoint to the present project is twofold. First, target-like states provide a natural comparison class for source-type behavior in one dimension: they are localized wave emitters with a structured core and outgoing waves in the far field \cite{Stich2003,KopellHoward1981}. Second, the pacemaker literature shows that wave-emitting states may arise through more than one mechanism, including imposed heterogeneity, self-organization, or bistability between distinct local dynamical states \cite{Stich2003}. For this reason, when studying numerical regimes in the one-dimensional Brusselator, it is natural to compare target-like states not only with more symmetric wave or stripe regimes, but also with spiral-like or source-defect-like states that organize the surrounding oscillatory medium in a different way \cite{Stich2003,SandstedeScheel2004}.

The second relevant strand of literature concerns one-dimensional ``spirals.'' Strictly speaking, true spiral waves are two-dimensional rotating structures, so the term must be interpreted carefully in one spatial dimension \cite{Barkley1994,Barkley1995,SandstedeScheel2021}. Classical geometrical theories of spiral waves make clear that genuine spiral organization is fundamentally tied to two-dimensional wave rotation and curvature effects, which is exactly why the one-dimensional setting considered here must be interpreted as an analogue rather than a literal spiral-wave problem \cite{Keener1986}. More broadly, the spiral-wave literature also shows that symmetry breaking is central to the emergence and organization of spiral-type dynamics, reinforcing the relevance of asymmetry-based diagnostics in the present one-dimensional setting \cite{LeBlanc2002}. Nevertheless, Perraud \emph{et al.}\ report the experimental observation of endogenous antisymmetric wave sources in a quasi-one-dimensional chemical system and explicitly interpret them as one-dimensional analogues of spiral organization \cite{PerraudBorckmans1993}. In their experiments and accompanying Brusselator simulations, these structures act as localized sources that emit waves asynchronously to the left and right, thereby creating a distinctly source-type or left--right asymmetric organization of the medium \cite{PerraudBorckmans1993}. This makes them highly relevant to the present work, where the spiral-like regime is likewise not a literal two-dimensional spiral, but rather a one-dimensional wave-emitting state with broken symmetry relative to more regular wave or stripe states \cite{PerraudBorckmans1993,CytrynbaumLewis2009}. Related studies of antispiral and wave-source behavior in oscillatory reaction--diffusion media further support the idea that wave-emitting structures can organize surrounding dynamics in ways that are not captured by the simplest symmetric wave-train picture \cite{NicolaBruschBaer2004}.

Equally important is the mechanism proposed by Perraud \emph{et al.}\ for the emergence of such structures. Their interpretation relies on the interaction between Turing and Hopf modes near a codimension-two setting, so that a localized Turing-like core may be embedded in an oscillatory background and act as a source of outward-traveling waves \cite{PerraudBorckmans1993}. This is especially suggestive for the Brusselator, where competition between oscillatory and spatially structured tendencies provides a natural mechanism for the appearance of mixed or defect-mediated states \cite{YuGumel2001,TzouMaBaylissMatkowskyVolpert2013}. In this way, the one-dimensional ``spirals'' of Perraud \emph{et al.}\ and the target/pacemaker viewpoint of Stich together motivate the two principal pattern classes considered later in this work \cite{PerraudBorckmans1993,Stich2003}.

Taken together, these literatures show that target-like and spiral-like or source-like states are not ad hoc numerical categories, but are well grounded in prior work on wave sources in oscillatory media \cite{Stich2003,PerraudBorckmans1993,SandstedeScheel2004}. Target patterns emphasize the role of localized pacemakers and outgoing waves \cite{Stich2003,KopellHoward1981}, while one-dimensional ``spirals'' show that asymmetric or source-like wave-emitting structures can arise even in settings without genuine two-dimensional spiral geometry \cite{PerraudBorckmans1993,CytrynbaumLewis2009}. These perspectives make it natural to study target-like and spiral/source-like regimes side by side, and they provide the conceptual background for the feature-based distinction developed in the later sections \cite{ZhaoMaffaSandstede2025}.

\begin{figure}[H]
    \centering
    \includegraphics[width=1.0\textwidth]{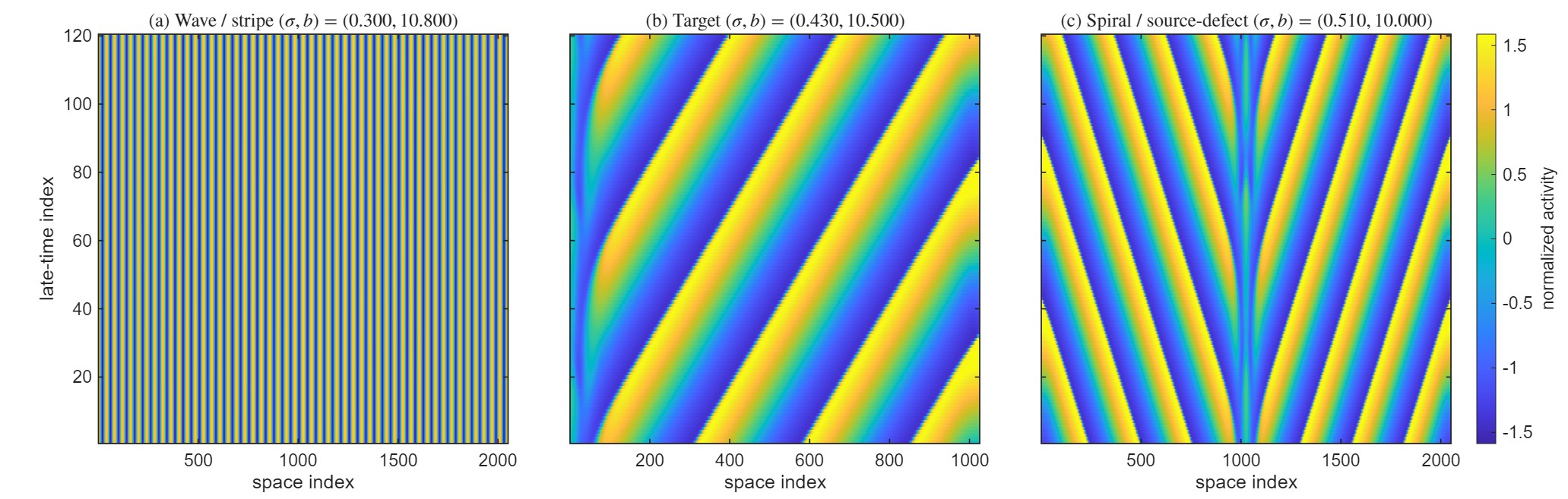}
    \caption{Representative late-time spacetime plots of the second component $v(t,x)$ in the one-dimensional Brusselator: (a) a wave/stripe-like state, (b) a target-like state (shown here in its half-domain representation), and (c) a spiral/source-defect-like state. These qualitative regimes motivate the transition problems studied in this work.}
    \label{fig:intro_regimes}
\end{figure}

Figure~\ref{fig:intro_regimes} shows representative late-time spacetime plots of the three qualitative regimes that motivate the rest of this work: a wave/stripe-like state, a target-like state, and a spiral/source-defect-like state.

In the one-dimensional setting, target and spiral/source-defect patterns can appear visually quite similar, especially when one plots only a half-domain target state \cite{PerraudBorckmans1993,Stich2003}. For this reason, the work distinguishes carefully between the geometric representation being shown and the underlying pattern class being discussed.

Throughout this manuscript, it is important to distinguish between the full target pattern, defined on a symmetric domain such as $[-L,L]$ with the core at $x=0$, and its half-domain representation on $[0,L]$. These are not different pattern types in any essential dynamical sense: the half-target is simply one side of the full target pattern, plotted on a reduced domain with the same core location \cite{Stich2003,KopellHoward1981}. Since the half-domain formulation is often more convenient for computation and feature design, both representations appear in this work. Whenever a target pattern is shown, we therefore indicate explicitly whether the plotted state is the full pattern or only its half-domain version.

\subsection{Goal of This Work}

This work studies a one-dimensional Brusselator model in the two-parameter plane $(\sigma,b)$, with the goal of understanding how qualitatively different late-time regimes are organized in parameter space. The main objects of interest are not isolated solution branches, but \emph{transition curves} separating regions in which different pattern types prevail. In particular, the focus is on transitions between spiral/source-defect-like and wave/stripe-like regimes, and between target-like and non-target-like regimes.

The computational viewpoint adopted here differs from classical continuation approaches for exact coherent structures \cite{BeynChampneysDoedelGovaertsKuznetsovSandstede2002,ChampneysSandstede2007,Uecker2021}. In standard settings, one formulates a boundary-value problem whose solutions represent the structure of interest and then applies Newton or pseudo-arclength continuation to trace solution branches or bifurcation curves \cite{ChampneysSandstede2007,Uecker2021}. By contrast, the present problem is formulated in a more data-driven way. Instead of solving directly for exact coherent structures, we evolve the underlying partial differential equation numerically, extract scalar features from the resulting late-time behavior, and use those features to detect regime changes \cite{ZhaoMaffaSandstede2025}. In this approach, a transition curve is represented as the threshold level set of a simulation-based scalar observable rather than as the zero set of an explicitly derived bifurcation equation \cite{ZhaoMaffaSandstede2025}.

The central methodological idea is therefore to treat suitably designed scalar observables as empirical order parameters \cite{ZhaoMaffaSandstede2025}. For the spiral transition, the relevant observable is a spacetime-based symmetry-defect feature that distinguishes asymmetric spiral/source-like behavior from more symmetric wave/stripe dynamics. For the target transition, the observable is a composite detector that identifies the simultaneous presence of a structured core and an oscillatory tail. Where these features give regular crossings, the associated level sets are tracked numerically. Where mixed patterns weaken this interpretation, the transition is estimated from direct parameter sweeps.

To map these transitions, this work develops a feature-based continuation framework built around a secant-based predictor and a local-sweep corrector \cite{PapakonstantinouTapia2013,AlvesDaSilvaCastroCosta2003,ZhaoMaffaSandstede2025}. The corrector performs a short local sweep in one parameter direction, locates the relevant threshold crossing, and selects the appropriate branch using crossing direction and hysteresis information. Because the spiral feature is branch-sensitive, branch-aware detector logic is incorporated into the algorithm. The continued curves are checked against direct spacetime diagnostics and fixed-parameter sweeps. In lower regions where the scalar features become less specific, vertical sweeps and direct pattern classification are used to estimate the transition instead of treating every threshold crossing as continuation output.

The goal of this work is therefore to develop and test a feature-based continuation framework for detecting and tracing transition curves between source/spiral-like and wave/stripe-like regimes, as well as between target-type and non-target-type regimes, in a one-dimensional Brusselator model \cite{ZhaoMaffaSandstede2025}.

\section{The Brusselator model and pattern regimes}

\subsection{The PDE Model}

In this work, we study a one-dimensional reaction--diffusion Brusselator model on a bounded interval \cite{Wazwaz2000,YuGumel2001,GolovinMatkowskyVolpert2008}. The unknowns are two scalar fields
\[
u=u(t,x), \qquad v=v(t,x),
\]
where $t \geq 0$ denotes time and $x \in [0,L]$ denotes the spatial variable. In the chemical interpretation of the Brusselator, $u$ and $v$ represent the two interacting species of the underlying reaction--diffusion system \cite{Murray2003,EpsteinPojman1998}. Throughout this work, the second component $v$ plays a particularly important role, since the principal feature functions introduced later are extracted from either its final spatial profile or its late-time spacetime field.

The model solved numerically is
\begin{equation}
\begin{aligned}
u_t &= d_1 u_{xx} + a - (b+1)u + u^2v,\\
v_t &= d_2 v_{xx} + bu - u^2v,
\end{aligned}
\label{eq:brusselator_1d}
\end{equation}
posed on the interval $[0,L]$ \cite{Wazwaz2000,YuGumel2001}. In the implementation used for this project, the constant
\[
a=2.5
\]
is fixed, as is the second diffusion coefficient
\[
d_2 = 9.73.
\]
The two-parameter plane of interest is $(\sigma,b)$, where $b$ is the usual Brusselator kinetic parameter and $\sigma$ is used to define the first diffusion coefficient through
\[
d_1 = \sigma d_2.
\]
Thus, the continuation problem is formulated in the parameter plane
\[
(\sigma,b),
\]
with all other constants held fixed.

The boundary conditions are homogeneous Neumann conditions,
\begin{equation}
u_x(t,0)=u_x(t,L)=0, \qquad v_x(t,0)=v_x(t,L)=0,
\label{eq:neumann_bc}
\end{equation}
which correspond to no-flux conditions at the endpoints of the interval \cite{Pao1982,Murray2003}. These boundary conditions are built directly into the finite-difference discretization by modifying the second-derivative operator at the first and last grid points. They are natural for the source-like and target-like wave-emitting states considered later, since they allow the interior dynamics to organize without imposing fixed Dirichlet values at the boundary \cite{PerraudBorckmans1993,Stich2003}.

Initial conditions are supplied in discretized form as
\[
u(0,x)=u_0(x), \qquad v(0,x)=v_0(x),
\]
with the specific choice depending on the regime under study. In some parameter regions, weakly perturbed initial data are sufficient to produce representative wave or stripe states \cite{Pearson1993,ZhaoMaffaSandstede2025}. In other regions, especially when one wishes to access target-like or spiral/source-like states, more structured seeds are used \cite{PerraudBorckmans1993,Stich2003,ZhaoMaffaSandstede2025}. The work therefore treats initial conditions pragmatically, as a way of generating representative late-time regimes whose boundaries in parameter space are the real objects of interest \cite{ZhaoMaffaSandstede2025}.

The numerical solution framework follows the method-of-lines approach \cite{Wazwaz2000,Uecker2021}. The spatial interval is discretized on a uniform grid
\[
x_1,\dots,x_N,
\]
with mesh width $\Delta x$, and the second derivative in \eqref{eq:brusselator_1d} is approximated by a finite-difference operator \cite{Wazwaz2000,Uecker2021}. In particular, a discrete Laplacian with Neumann boundary modifications is used, producing a coupled system of ordinary differential equations in time for the nodal values of $u$ and $v$ \cite{Pao1982,Uecker2021}. This semi-discrete system is then integrated forward in time using \texttt{ode45} in MATLAB. In the solver used here, the solution is recorded on a prescribed time grid with spacing $0.1$, and the second component is stored as a spacetime matrix
\[
V = \bigl(v(t_j,x_i)\bigr)_{j,i},
\]
which is subsequently used for visualization, diagnostics, and feature extraction.

At a high level, the role of the PDE model in this work is not only to define the underlying reaction--diffusion dynamics, but also to generate the late-time data from which scalar observables are built \cite{ZhaoMaffaSandstede2025}. The transition curves studied later are therefore not obtained from an explicit analytical bifurcation equation. Instead, they are inferred from how the long-time behavior of \eqref{eq:brusselator_1d} changes as $(\sigma,b)$ varies \cite{ZhaoMaffaSandstede2025}. For this reason, the one-dimensional Brusselator serves here both as the dynamical system of interest and as the simulation engine underlying the feature-based continuation framework \cite{ZhaoMaffaSandstede2025}.

\subsection{Pattern Regimes of Interest}

With the model now specified, we turn to the principal qualitative regimes that appear in direct numerical simulations of \eqref{eq:brusselator_1d}. As discussed in Section~1, the present work is concerned with wave/stripe-like, spiral/source-like, and target-like behavior. The purpose of this section is therefore not to reintroduce these regimes from a broad pattern-formation perspective, but rather to describe how they appear \emph{numerically} in the one-dimensional Brusselator and how they are interpreted in the later continuation framework.

\begin{figure}[H]
    \centering
    \includegraphics[width=0.7\textwidth]{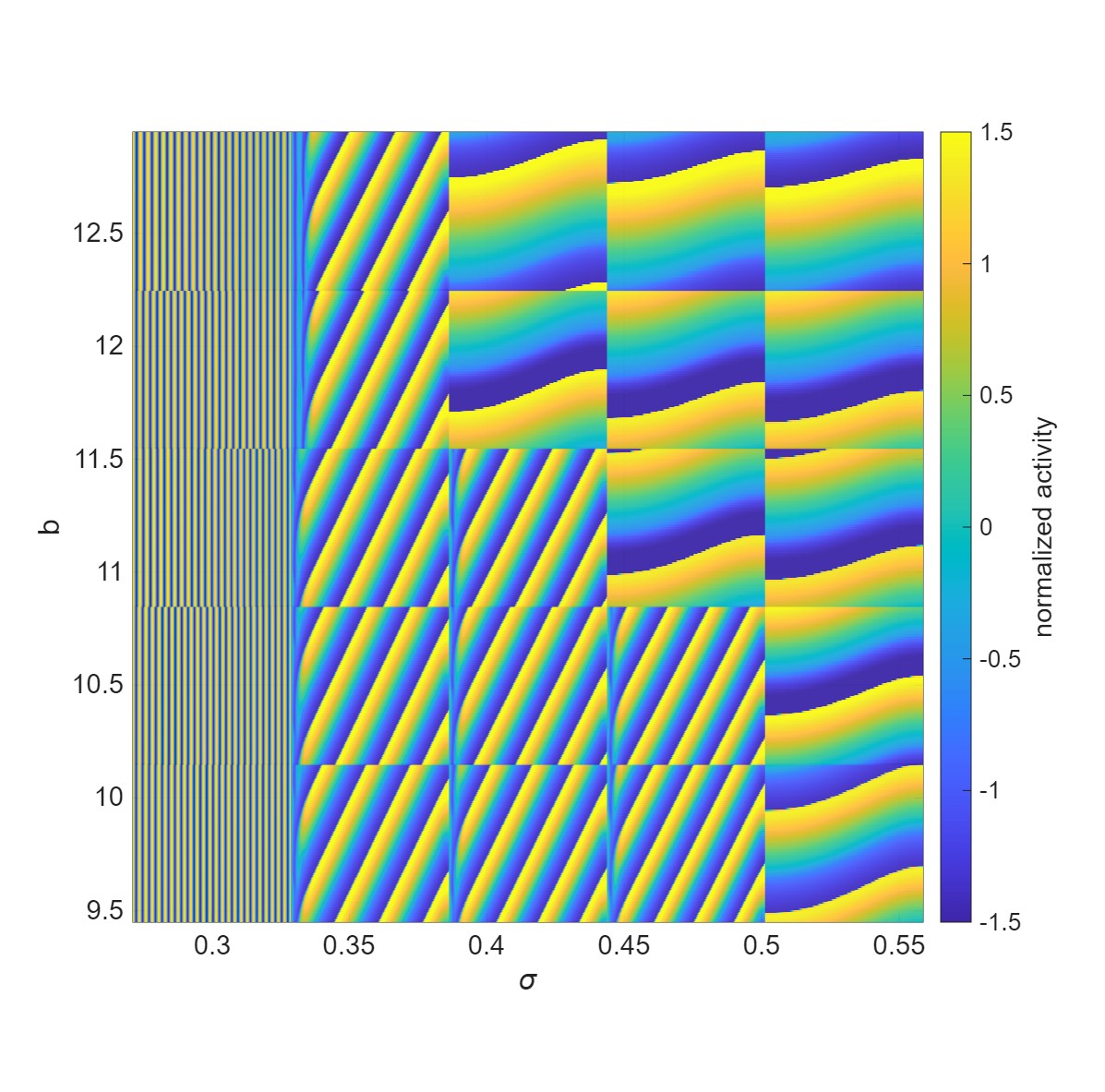}
    \caption{Qualitative late-time regime atlas for representative parameter values in the $(\sigma,b)$-plane. Each tile shows the late-time spacetime field of the second component for a sampled parameter pair. The atlas illustrates the coexistence of several qualitatively distinct regimes, including wave/stripe-like behavior, target-like states, spiral/source-defect-like states, and intermediate patterns.}
    \label{fig:regime_atlas}
\end{figure}

A key point is that the regimes considered here are \emph{late-time numerical regimes}, not exact analytically characterized solution classes \cite{ZhaoMaffaSandstede2025,Stich2003}. In other words, the labels ``wave/stripe-like,'' ``spiral/source-like,'' and ``target-like'' refer to robust qualitative behaviors observed in long-time simulations of the PDE, typically through spacetime plots of the second component $v(t,x)$ and related diagnostics \cite{ZhaoMaffaSandstede2025}. This distinction is important because the transition curves studied later are not defined by closed-form bifurcation formulae or exact coherent-structure equations, but by empirical changes in these observed late-time behaviors \cite{ZhaoMaffaSandstede2025}. To give a coarse global view of these late-time numerical regimes, Figure~\ref{fig:regime_atlas} shows a qualitative atlas of representative spacetime patterns across a sampled subset of the $(\sigma,b)$-plane.

Figure~\ref{fig:regime_atlas} serves as a qualitative overview showing that the one-dimensional Brusselator supports several visibly distinct late-time regimes across the sampled parameter range, and that these regimes are organized nontrivially in the $(\sigma,b)$-plane \cite{ZhaoMaffaSandstede2025}.

The most regular regime encountered in the simulations is the \emph{wave/stripe-like} regime. In spacetime plots, this regime typically appears as a comparatively symmetric and spatially organized pattern, often resembling traveling or weakly modulated wave trains \cite{SandstedeScheel2004,DoelmanSandstedeScheelSchneider2009}. Depending on parameter values and the precise seed, the pattern may look more stripe-like or more wave-like, but in either case the behavior is comparatively regular and lacks the clearly localized wave-emitting core that characterizes the other two regimes \cite{ZhaoMaffaSandstede2025}. This regime serves as the natural comparison class for the spiral/source-like transition studied later.

A second important regime is the \emph{spiral/source-like} regime. Here the late-time spacetime plot exhibits a localized organizing structure that emits or mediates waves asymmetrically, producing a source-defect-like appearance rather than a fully symmetric wave train \cite{PerraudBorckmans1993,SandstedeScheel2004}. In the one-dimensional setting, this is not a literal two-dimensional spiral, but rather a source-like or defect-like state whose spacetime organization breaks the more regular symmetry seen in the wave/stripe regime \cite{PerraudBorckmans1993,CytrynbaumLewis2009}. It is precisely this qualitative distinction that motivates the spacetime symmetry-defect feature introduced in Section~4.

The third regime of interest is the \emph{target-like} regime. Numerically, this regime is characterized by a localized core together with outgoing wave activity in the surrounding medium \cite{KopellHoward1981,Stich2003}. In the computations used in this work, target-like states are often associated with half-target or source-defect seeds and are recognized not through exact geometric circularity, which would be meaningless in one dimension, but through the simultaneous presence of a structured core and an oscillatory tail \cite{Stich2003}. This motivates the composite target feature introduced later in Section~5.

In addition to these three principal classes, one also encounters \emph{mixed} or \emph{ambiguous} regimes in certain parts of parameter space, especially near transition boundaries or in parameter regions where the organizing structure is weak \cite{ZhaoMaffaSandstede2025}. These states may exhibit partial source-like behavior without a clean localized core, or may display irregular mixtures of wave-train and defect-like dynamics \cite{PerraudBorckmans1993,ZhaoMaffaSandstede2025}. Such regimes are important because they explain why neither visual inspection nor a scalar detector alone is sufficient in every part of parameter space \cite{ZhaoMaffaSandstede2025}. They motivate the use of carefully designed scalar observables in regular regions together with direct spacetime classification where the regime distinction becomes ambiguous \cite{ZhaoMaffaSandstede2025}.

For later reference, Table~\ref{tab:regime_summary} summarizes the main numerical characteristics of the pattern classes considered in this work.

\begin{table}[H]
\centering
\caption{Qualitative summary of the principal late-time numerical regimes observed in the one-dimensional Brusselator.}
\label{tab:regime_summary}
\begin{tabular}{p{3.2cm} p{4.2cm} p{3.8cm} p{3.6cm}}
\toprule
\textbf{Regime} & \textbf{Late-time spacetime appearance} & \textbf{Localized organizing core} & \textbf{Role in this work} \\
\midrule
Wave / stripe-like
&
Comparatively regular, symmetric, wave-train or stripe-like spacetime organization
&
Typically absent
&
Reference regime for spiral/source-like transition \\[0.6em]

Spiral / source-like
&
Asymmetric spacetime organization with source-defect-like wave emission
&
Present in defect/source form
&
Detected by branch-aware spacetime symmetry-defect feature \\[0.6em]

Target-like
&
Localized core together with outgoing oscillatory tail behavior
&
Present
&
Detected by composite core/tail feature \\[0.6em]

Mixed / ambiguous
&
Irregular or partially organized spacetime behavior near regime boundaries
&
Weak, unclear, or unstable
&
Explains the need for feature-based rather than purely visual transition detection \\
\bottomrule
\end{tabular}
\end{table}

Representative spacetime plots for these regimes will be used throughout the work, beginning with the introductory figure in Section~1 and continuing in the numerical results section. At this stage, the main point is that the one-dimensional Brusselator supports several qualitatively distinct late-time behaviors in the $(\sigma,b)$ plane, and that these behaviors are sufficiently robust to motivate the search for transition curves separating their regions of prevalence \cite{ZhaoMaffaSandstede2025,PerraudBorckmans1993,Stich2003}. This leads naturally to the next section, where the emphasis shifts from the regimes themselves to the question of why their boundaries in parameter space are the primary objects of study.

\subsection{Why Transition Curves Matter}

The central object of study in this work is not an isolated simulation at a single parameter value, but the organization of qualitative dynamics across the two-parameter plane $(\sigma,b)$ \cite{CrossHohenberg1993,ZhaoMaffaSandstede2025}. As the parameters vary, the one-dimensional Brusselator exhibits different prevailing late-time regimes, including wave/stripe-like, spiral/source-like, and target-like behavior. The natural scientific question is therefore not only \emph{what} patterns occur, but also \emph{where} in parameter space one regime gives way to another \cite{CrossHohenberg1993,ZhaoMaffaSandstede2025}. In other words, one seeks to determine the boundaries separating regions with qualitatively different long-time dynamics.

This viewpoint is well motivated by the broader literature on pattern formation. In reaction--diffusion systems and other nonequilibrium media, different parameter regions may support different prevalent patterns, and the geometry of these regions encodes important dynamical information \cite{CrossHohenberg1993,Hoyle2006,ZhaoMaffaSandstede2025}. Zhao, Maffa, and Sandstede formulate this issue explicitly: different patterns may prevail in different parts of parameter space, and the boundaries between these regions correspond to transition or bifurcation curves that organize the global pattern-selection picture \cite{ZhaoMaffaSandstede2025}. From this perspective, a coarse collection of example simulations is only a first step; what one ultimately wants is a more systematic description of the regime boundaries themselves \cite{ZhaoMaffaSandstede2025}.

A natural first step is to simulate the PDE on a parameter grid and inspect the resulting patterns visually \cite{CrossHohenberg1993,ZhaoMaffaSandstede2025}. This gives a useful coarse overview, but it becomes unreliable near narrow or curved transition regions and in areas where patterns look similar or mixed states occur \cite{ZhaoMaffaSandstede2025}. As Zhao \emph{et al.}\ emphasize, direct parameter scans are valuable for a qualitative picture but limited in accuracy near regime boundaries \cite{ZhaoMaffaSandstede2025}.

Classical continuation offers a more precise alternative when the relevant states can be formulated as a boundary-value problem \cite{BeynChampneysDoedelGovaertsKuznetsovSandstede2002,ChampneysSandstede2007,Uecker2021}. In that setting, one can compute equilibria, waves, or bifurcation curves directly \cite{ChampneysSandstede2007,Uecker2021}. Here, however, the transitions are identified through late-time simulations and pattern diagnostics rather than through an explicit analytical characterization of a coherent structure \cite{ZhaoMaffaSandstede2025}. As Zhao \emph{et al.}\ note, this makes classical continuation harder to apply directly when the goal is to trace dynamically selected pattern transitions \cite{ZhaoMaffaSandstede2025}.

Rather than deriving an exact analytical bifurcation equation for each transition, we assign to each simulation a scalar observable that captures a regime-defining property of the late-time solution \cite{ZhaoMaffaSandstede2025}. We then define the transition set in the $(\sigma,b)$-plane as a threshold level set of this observable. The resulting boundaries are therefore \emph{feature-defined transition curves}, not classical analytical bifurcation curves \cite{ZhaoMaffaSandstede2025}. Their interpretation is empirical but still mathematically structured, since they arise as level sets of scalar quantities extracted systematically from the PDE dynamics \cite{ZhaoMaffaSandstede2025}.

This distinction is important. A feature-defined transition curve is meant to separate regions with different prevailing late-time behavior, even when no explicit instability condition or exact coherent-structure equation is available \cite{ZhaoMaffaSandstede2025}. In this work, the feature functions introduced later serve as numerical order parameters: one distinguishes spiral/source-like behavior from wave/stripe-like behavior, and another detects target-like states. The continuation problem is then to trace threshold crossings of these observables across parameter space \cite{ZhaoMaffaSandstede2025}.

\begin{figure}[H]
    \centering
    \includegraphics[width=1\textwidth]{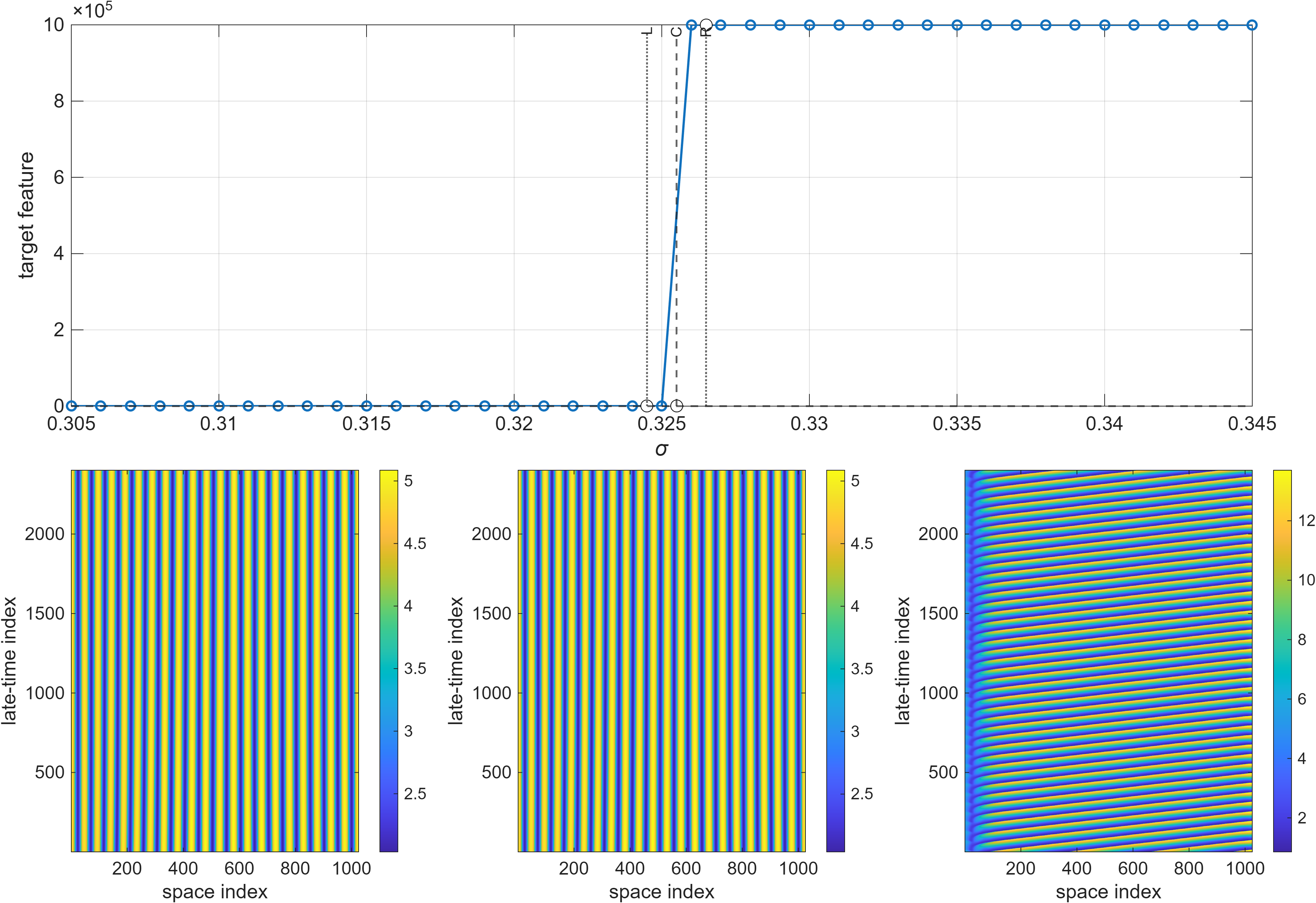}
    \caption{Motivating example for feature-defined transition curves. Panel (a) shows the target feature along a horizontal sweep at fixed $b=10.02$. Panels (b)--(d) show representative late-time spacetime plots to the left, near, and to the right of the detected transition. The sharp change in the scalar feature corresponds to a genuine change in qualitative late-time behavior, illustrating why it is meaningful to define and continue transition curves in parameter space.}
    \label{fig:target_sweep_motivation}
\end{figure}

A representative example is shown in Figure~\ref{fig:target_sweep_motivation}, where a horizontal sweep at fixed $b$ produces a sharp change in the target feature that corresponds to a real qualitative change in the late-time spacetime pattern. This is precisely why it is useful to regard the transition set as a geometric object in parameter space and to compute it systematically rather than relying only on isolated simulations or coarse parameter scans.

In summary, transition curves matter because they provide a compact and geometrically meaningful description of how different late-time regimes are arranged in the $(\sigma,b)$-plane \cite{CrossHohenberg1993,ZhaoMaffaSandstede2025}. Direct parameter scans alone are too coarse for this purpose \cite{ZhaoMaffaSandstede2025}, while classical bifurcation methods are not always directly applicable to the transitions considered here \cite{ChampneysSandstede2007,Uecker2021,ZhaoMaffaSandstede2025}. The approach taken in this work is therefore to compute transition curves numerically through scalar feature thresholds, thereby obtaining a data-driven description of regime boundaries in the one-dimensional Brusselator \cite{ZhaoMaffaSandstede2025}.

\section{Continuation Framework}

\subsection{Transition Curves as Feature-Defined Level Sets}

The continuation problem considered in this work begins with a simple geometric idea. Rather than computing a single distinguished solution at one parameter value, we seek a \emph{curve} in the parameter plane
\[
p=(\sigma,b)\in\mathbb{R}^2
\]
that separates two qualitatively different late-time regimes of the PDE \cite{CrossHohenberg1993,ZhaoMaffaSandstede2025}. On one side of this curve, simulations converge to one type of behavior; on the other side, to another. The transition curve is therefore the set of parameter values at which a chosen transition condition holds.

In classical continuation theory, such a curve is often defined implicitly by an equation
\[
F(p)=0,
\]
where $F$ is an analytically defined bifurcation or solvability condition \cite{BeynChampneysDoedelGovaertsKuznetsovSandstede2002,ChampneysSandstede2007,Uecker2021}. Here, however, the relevant quantity is not available in closed form. Instead, it is constructed from direct numerical simulations. This viewpoint is closely aligned with the data-driven continuation framework of Zhao, Maffa, and Sandstede, where boundaries between regions with different prevailing patterns are inferred from simulation-based observables rather than from an explicit analytical equation \cite{ZhaoMaffaSandstede2025}.

Accordingly, the basic object in this work is a scalar feature map obtained from the PDE dynamics,
\[
F(\sigma,b),
\]
or, in threshold form,
\[
r(\sigma,b).
\]
Its meaning depends on the pattern class under consideration. For the spiral transition, it is a spacetime symmetry-defect feature; for the target transition, it is a composite observable built from core and tail variance measurements. In either case, the feature is produced by a simulation pipeline of the form
\[
(\sigma,b)
\longmapsto
\text{PDE solution}
\longmapsto
\text{late-time data}
\longmapsto
\text{scalar feature},
\]
so the continuation framework is built around an \emph{empirical response map} in parameter space \cite{ZhaoMaffaSandstede2025}.

If the transition is described by the zero set of a scalar feature, then the transition curve may be written abstractly as
\[
\Gamma=\{(\sigma,b):F(\sigma,b)=0\}.
\]
In the threshold formulation used here, one defines a scalar observable $r(\sigma,b)$ and a threshold value $r_{\mathrm{thr}}$, and sets
\[
\Gamma
=
\{(\sigma,b):r(\sigma,b)=r_{\mathrm{thr}}\}.
\]
Equivalently, one may define
\[
F(\sigma,b)=r(\sigma,b)-r_{\mathrm{thr}},
\]
so that the transition curve is again the zero level set of a scalar function.

This is an \emph{empirical level-set formulation} \cite{ZhaoMaffaSandstede2025}. The quantity $F$ is not an exact analytical bifurcation equation, but a simulation-derived scalar observable whose threshold is chosen to separate regimes robustly. For this reason, the curves computed here should be interpreted as \emph{feature-defined transition curves}, not exact analytical bifurcation curves \cite{CrossHohenberg1993,ZhaoMaffaSandstede2025}.

Once the transition is formulated as a level set in the $(\sigma,b)$-plane, the next question is how to trace it numerically. A naive approach would be to represent the transition as a graph over one parameter, for example
\[
b=b(\sigma)
\qquad\text{or}\qquad
\sigma=\sigma(b).
\]
However, this is generally too restrictive. Even if the curve is smooth, it may contain folds or turning points, in which case it ceases locally to be single-valued in one coordinate \cite{AllgowerGeorg2000,BeynChampneysDoedelGovaertsKuznetsovSandstede2002,Uecker2021}. Near such points, graph-based continuation is poorly aligned with the geometry of the transition set \cite{AllgowerGeorg2000,ChampneysSandstede2007}.

This issue is especially relevant here, since the continued object is a threshold level set of a simulation-based feature rather than an exact solution branch. There is no reason to expect such a set to be globally representable as a graph over $\sigma$ or $b$. Instead, the correct geometric viewpoint is to treat the transition as a parametrized path
\[
p(s)=(\sigma(s),b(s)),
\]
where $s$ is an arclength parameter, or more generally an approximate path-length parameter \cite{BeynChampneysDoedelGovaertsKuznetsovSandstede2002,ChampneysSandstede2007,Uecker2021}.

This viewpoint motivates the continuation framework developed below. Rather than stepping in only one coordinate direction and solving for the other, the algorithm uses previously accepted points to predict a new point in the local tangent direction, and then corrects that prediction by locating a nearby threshold crossing \cite{PapakonstantinouTapia2013,AlvesDaSilvaCastroCosta2003}. In this respect, the method is inspired by arclength continuation, but adapted to a setting where the underlying scalar function is obtained through simulation and feature extraction \cite{BeynChampneysDoedelGovaertsKuznetsovSandstede2002,ChampneysSandstede2007,ZhaoMaffaSandstede2025}.

\subsection{Secant Predictor and Local Sweep Corrector}

Once the transition set is formulated as a feature-defined level set in parameter space, the next step is to trace it numerically. The method used here is a predictor--corrector scheme, conceptually similar to arclength continuation: previously computed points predict a new point along the local tangent direction, and a correction step returns that point to the transition set \cite{BeynChampneysDoedelGovaertsKuznetsovSandstede2002,ChampneysSandstede2007,Uecker2021}. The key difference from classical continuation is that the implementation is not based on Newton solves for a boundary-value problem.

In the classical setting, coherent structures are continued as regular zeros of a suitably posed nonlinear system, often with phase conditions to remove symmetry-related degeneracies \cite{SandstedeNotes,ChampneysSandstede2007}. Here, by contrast, the continued object is a threshold level set of a scalar feature extracted from time-dependent simulations \cite{ZhaoMaffaSandstede2025}. Thus, the predictor--corrector philosophy is retained, but the corrector must be implemented differently.

Let
\[
p_k = (\sigma_k,b_k)
\]
denote an accepted point on the transition curve, and suppose that the previous point
\[
p_{k-1} = (\sigma_{k-1},b_{k-1})
\]
is also known. A natural approximation of the local tangent direction is the secant vector
\[
\tau_k
=
\frac{p_k-p_{k-1}}{\|p_k-p_{k-1}\|}.
\]
Given a step size $\Delta s >0$, the predictor step is
\[
p_{k+1}^{\mathrm{pred}}
=
p_k + \Delta s\, \tau_k.
\]
This is simply the rule ``continue in approximately the same direction'' \cite{PapakonstantinouTapia2013,AlvesDaSilvaCastroCosta2003}. Its role is geometric rather than analytic.

Write the predicted point as
\[
p_{k+1}^{\mathrm{pred}} = (\sigma^{\mathrm{pred}}, b^{\mathrm{pred}}).
\]
In a full pseudo-arclength method, one would now solve an augmented nonlinear system involving both the defining equation for the curve and an additional hyperplane condition near the predicted point \cite{BeynChampneysDoedelGovaertsKuznetsovSandstede2002,ChampneysSandstede2007}. That is not what is done here. Instead, the corrector takes the form of a \emph{local one-dimensional parameter sweep}: one fixes
\[
b = b^{\mathrm{pred}}
\]
and searches locally in the $\sigma$-direction near $\sigma^{\mathrm{pred}}$.

More precisely, one samples
\[
\sigma \in [\sigma^{\mathrm{pred}}-h,\;\sigma^{\mathrm{pred}}+h]
\]
for a chosen sweep half-width $h>0$. At each sampled value of $\sigma$, one runs the PDE simulation, extracts the relevant late-time feature, and evaluates
\[
r(\sigma,b^{\mathrm{pred}})
\]
or equivalently
\[
g_k(\sigma) := r(\sigma,b^{\mathrm{pred}})-r_{\mathrm{thr}}.
\]
The task of the corrector is then to locate a nearby threshold crossing
\[
g_k(\sigma_\ast)=0.
\]
Once such a crossing is identified, the corrected point is taken to be
\[
p_{k+1} = (\sigma_\ast,b^{\mathrm{pred}}).
\]

This local sweep corrector is the essential adaptation that makes continuation possible in the present data-driven setting \cite{ZhaoMaffaSandstede2025}. Instead of solving a two-dimensional nonlinear system by Newton iteration, the algorithm approximates the local intersection of the transition curve with the horizontal line
\[
b=b^{\mathrm{pred}}.
\]
The corrector therefore requires only that the scalar feature can be evaluated from simulation at nearby parameter values \cite{ZhaoMaffaSandstede2025}.

At the same time, the method differs from classical continuation. In the boundary-value setting, the object of interest is corrected by solving a regular nonlinear system directly \cite{ChampneysSandstede2007,Uecker2021}. Here, the corrector is one-dimensional and simulation-based: it fixes one parameter, samples the feature locally in the other, and selects the relevant threshold crossing \cite{ZhaoMaffaSandstede2025}. For this reason, the method is best described as a \emph{feature-based predictor--corrector continuation method inspired by arclength continuation} rather than a full pseudo-arclength Newton scheme.

\begin{figure}[H]
    \centering
    \includegraphics[width=0.72\textwidth]{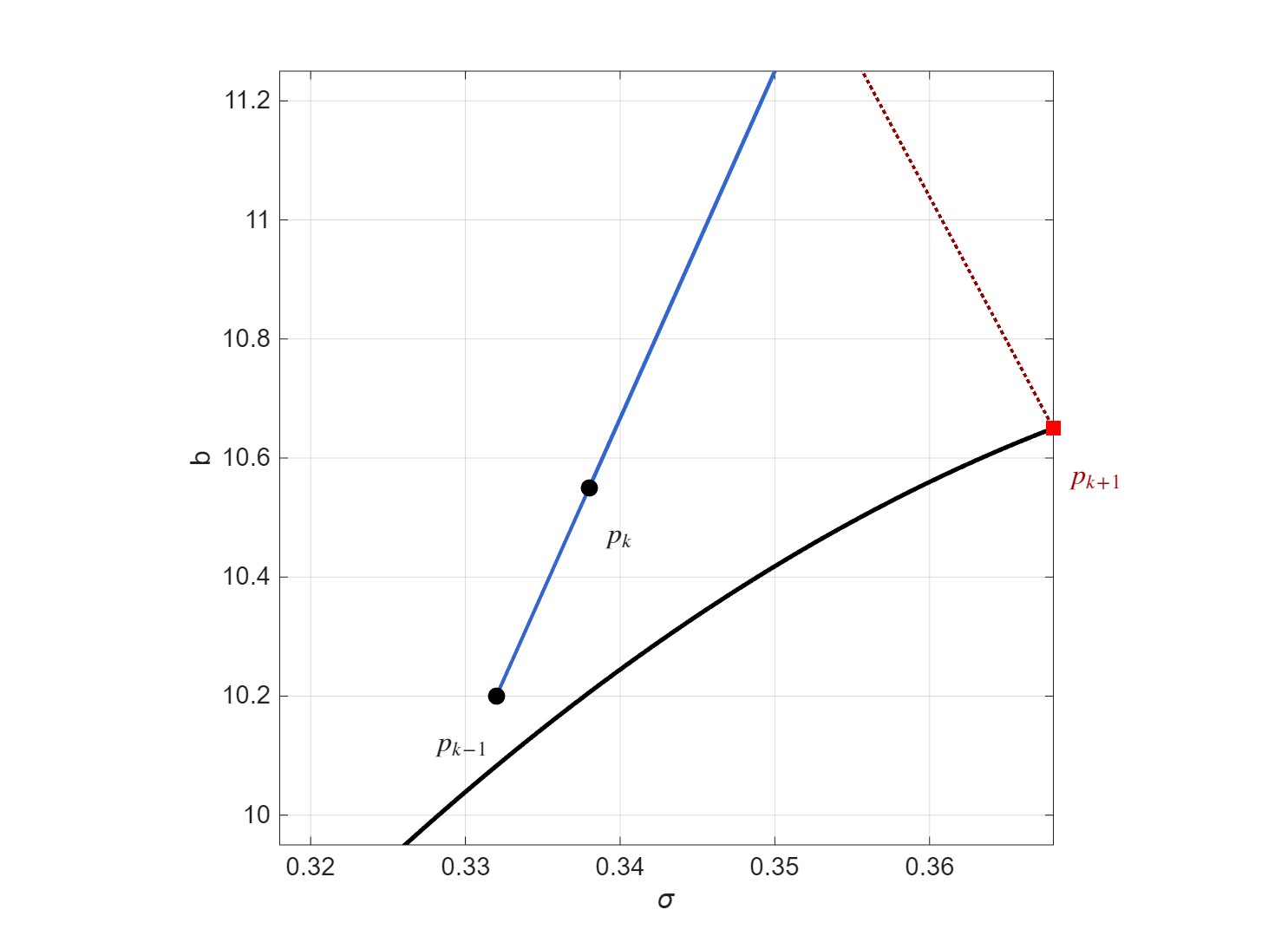}
    \caption{Schematic illustration of the secant predictor and local horizontal corrector. Given two previously accepted points $p_{k-1}$ and $p_k$ on the transition curve, the secant direction provides a local tangent approximation, and the next point is predicted by stepping forward along that direction. The black curve represents the transition set in the $(\sigma,b)$-plane.}
    \label{fig:secant_predictor_schematic}
\end{figure}

Figure~\ref{fig:secant_predictor_schematic} illustrates the geometric idea of the predictor step. The secant predictor provides a local estimate of where the next point on the branch should lie, without forcing the transition curve to be represented as a graph over a single coordinate \cite{BeynChampneysDoedelGovaertsKuznetsovSandstede2002,ChampneysSandstede2007}. The predictor advances tangentially along the branch, while the corrector projects that tentative point back onto the feature-defined transition set by a local threshold search \cite{ZhaoMaffaSandstede2025}. This construction is flexible enough to follow curves that bend or fold in parameter space, provided the local sweep intersects the transition curve in a sufficiently regular way \cite{AllgowerGeorg2000,Uecker2021}.

\subsection{Branch Selection, Hysteresis, Warm-Start versus Cold-Start, and Geometric Refinements}

The local sweep corrector described above does not merely need to find \emph{a} threshold crossing; it must find the \emph{correct} one associated with the branch being continued. In practice, this requires several refinements: branch-aware crossing direction, hysteresis to prevent jumping between nearby crossings, a careful distinction between warm-start and cold-start, and, in the lower part of parameter space, modifications to the continuation geometry itself \cite{ZhaoMaffaSandstede2025}.

A first issue is that the scalar feature may cross its threshold in different directions on different branches. Along one branch, the feature may decrease through the threshold as $\sigma$ increases, while along another branch it may increase through it. If one were to accept any threshold crossing indiscriminately, the corrector could lock onto the wrong side of the transition set. For this reason, the continuation algorithm imposes a \emph{directional crossing rule}. In the side-branch geometry used for most of the spiral and target computations, the selected crossing must satisfy not only
\[
r(\sigma,b)=r_{\mathrm{thr}},
\]
but also a local sign condition on the derivative with respect to the sweep variable, at least heuristically. For example, one branch may require a downward crossing with
\[
\frac{\partial r}{\partial \sigma}<0,
\]
while another requires an upward crossing with
\[
\frac{\partial r}{\partial \sigma}>0.
\]
Thus, branch selection is encoded not only by the threshold itself, but also by the direction in which the threshold is crossed.

A second issue arises when the local sweep contains multiple admissible threshold crossings. This can happen if the feature oscillates, develops spikes, or if more than one regime transition is present in the sweep window. In that situation, selecting a crossing solely by threshold equality is unstable: the corrector may jump between nearby crossings even when the underlying branch should vary smoothly. To reduce this ambiguity, the algorithm incorporates a \emph{hysteresis rule}. If the previously accepted point on the branch has corrected coordinate $\sigma_k$, and the current sweep produces candidate crossings
\[
\sigma^{(1)},\dots,\sigma^{(m)},
\]
then the algorithm prefers the crossing nearest to the previous accepted location, namely the one minimizing
\[
\bigl|\sigma^{(j)}-\sigma_k\bigr|.
\]
This is not a classical analytical condition, but a practical regularization that enforces continuity of the numerically detected branch.

A third issue concerns dependence on the initial condition used for the PDE simulation. Since the scalar observable is obtained by simulation rather than given explicitly, it is more accurate to think of it as
\[
r(\sigma,b;u_0),
\]
where $u_0$ denotes the seed used to generate the late-time state. This makes the distinction between warm-start and cold-start essential. Along an accepted continuation branch, warm-start is desirable: once a valid transition point has been found, the final state at that point provides a natural seed for nearby simulations, helping the computation remain on the same solution family. By contrast, using warm-start \emph{inside the local sweep} changes the mathematical meaning of the sampled feature. If the simulation at one value of $\sigma$ is initialized from the final state at the previous sampled value, then the resulting map is no longer a consistent function
\[
\sigma \longmapsto r(\sigma,b^{\mathrm{pred}}),
\]
but rather a path-dependent object. For this reason, the local corrector uses a fixed detection seed, or cold-start, so that the sweep more faithfully samples a well-defined response curve. Thus,
\begin{itemize}
    \item warm-start is useful \emph{along} the continuation branch,
    \item cold-start is preferable \emph{within} the local detection sweep.
\end{itemize}

A further issue arose in the lower part of parameter space, where the computed left and right downward branches appeared to narrow toward one another. This suggested that the transition region might not consist only of left and right side walls, but could also contain a lower connecting piece. The side branches were naturally treated by a horizontal corrector at fixed $b$, since locally they behaved like
\[
\sigma=\sigma(b).
\]
The suspected lower connecting piece is more naturally viewed as
\[
b=b(\sigma),
\]
for which a horizontal corrector is poorly aligned with the geometry. This observation motivated vertical sweeps in $b$ at selected fixed values of $\sigma$. At each selected value of $\sigma$, the simulations were started from the same prescribed seed, and the late-time spacetime plots were inspected as $b$ was varied. The transition was placed between neighboring sampled values of $b$ for which the observed pattern changed class. The middle points shown later are based on these visually identified transition brackets. Lines drawn between them are interpolation guides and are not additional continuation output.

A vertical-corrector continuation was also tested in this region by fixing $\sigma$ and searching for a threshold crossing in $b$. This test showed why the distinction is important. In mixed or irregular regions, the scalar detector may cross its nominal threshold even when the corresponding pattern change is not the intended spiral-to-wave or target-to-wave transition. A smooth threshold curve can therefore continue into a different pattern class. For this reason, the vertical-corrector output is treated as a diagnostic of the feature rather than as the source of the lower boundary reported in the regime diagrams.

Temporal coherence provides an additional diagnostic for this ambiguity. Let $V$ denote the recorded spacetime matrix of the second component, and restrict to the last portion $V_{\mathrm{late}}$ of the observation window. Define adjacent-frame differences by
\[
d_k
=
\frac{\|V_{\mathrm{late}}(k,:)-V_{\mathrm{late}}(k-1,:)\|}
{\|V_{\mathrm{late}}(k,:)\|+\varepsilon},
\]
and let
\[
r_{\mathrm{coh}}
=
\frac{1}{M-1}\sum_{k=2}^{M} d_k
\]
be their mean over the retained time window. The pattern is regarded as temporally coherent if
\[
r_{\mathrm{coh}} < r_{\mathrm{coh}}^{\mathrm{thr}}.
\]
Values that fail this check are flagged for direct inspection rather than accepted from the scalar threshold alone. The coherence score does not by itself classify a pattern, but it helps identify parameter values at which a visually based transition bracket is more appropriate than automatic continuation.

\begin{figure}[H]
    \centering
    \includegraphics[width=0.82\textwidth]{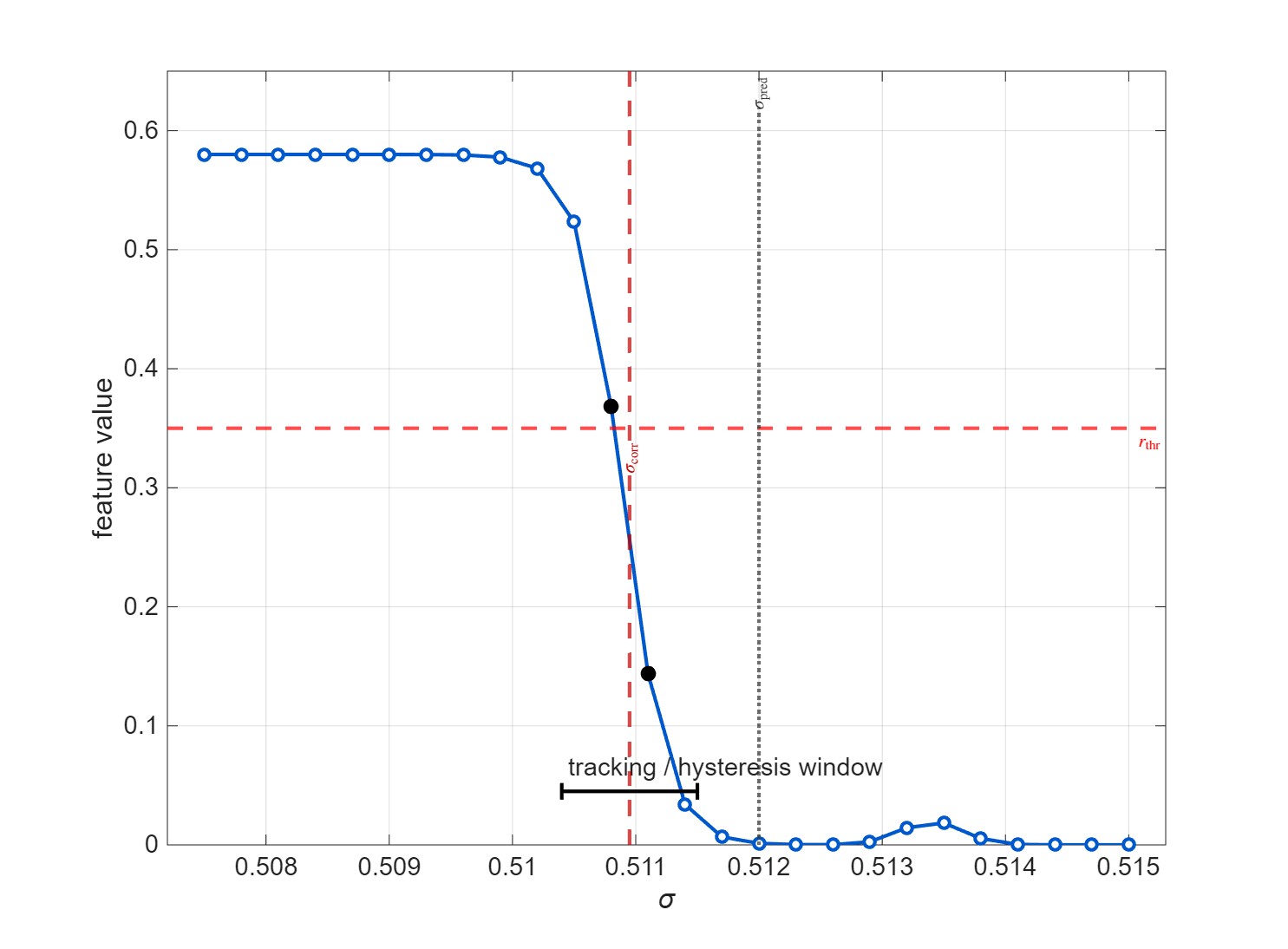}
    \caption{Schematic illustration of local sweep correction and threshold-crossing selection. The horizontal dashed line indicates the feature threshold $r_{\mathrm{thr}}$, while the vertical dotted line marks the predicted value $\sigma^{\mathrm{pred}}$. The corrected point is chosen from a nearby threshold crossing, with the tracking/hysteresis window used to prefer continuity with the previously accepted branch.}
    \label{fig:corrector_hysteresis_schematic}
\end{figure}

The threshold-selection logic used by the local corrector is illustrated schematically in Figure~\ref{fig:corrector_hysteresis_schematic}. In practice, a local sweep may contain more than one admissible crossing, and without additional selection logic the corrector may jump from one crossing to another. The hysteresis rule regularizes this process by favoring the crossing most consistent with the previously accepted branch location.

Taken together, these refinements show that the continuation problem in this work is not merely a matter of following a numerical level set. It is a branch-tracking problem in which the scalar feature, the corrector geometry, the crossing direction, and the regularity of the underlying pattern all interact \cite{ChampneysSandstede2007,Uecker2021,ZhaoMaffaSandstede2025}. The resulting framework is therefore more structured than a raw threshold search, but also more flexible than a classical boundary-value continuation method \cite{ZhaoMaffaSandstede2025}.

\subsection{Relation to Classical Pseudo-Arclength Continuation}

It is useful to state clearly how the continuation method used in this work relates to classical pseudo-arclength continuation \cite{BeynChampneysDoedelGovaertsKuznetsovSandstede2002,ChampneysSandstede2007,Uecker2021}. In the standard setting, one seeks zeros of an explicitly formulated nonlinear system and applies a predictor--corrector scheme in which the corrector solves an augmented system consisting of the defining equation together with a hyperplane condition near the predicted point \cite{BeynChampneysDoedelGovaertsKuznetsovSandstede2002,ChampneysSandstede2007}. This framework is especially effective when the object of interest can be written as a coherent structure satisfying a well-posed boundary-value problem \cite{SandstedeNotes,ChampneysSandstede2007}. Champneys and Sandstede emphasize that, in such cases, phase conditions and Newton-based correction provide a robust way to continue solution branches and organize them in parameter space \cite{ChampneysSandstede2007}. From the broader viewpoint of bifurcation theory, this classical setting belongs to a well-developed framework in which branches are continued as solution sets of nonlinear equations posed in function spaces \cite{ChowHale1982}.

The method used here retains the geometric predictor--corrector idea, but not the full classical corrector. The predictor is secant-based, as in arclength continuation \cite{PapakonstantinouTapia2013,AlvesDaSilvaCastroCosta2003}, but the correction step does not solve an augmented two-dimensional nonlinear system. Instead, it performs a one-dimensional local sweep at fixed $\sigma$ or fixed $b$, evaluates a simulation-derived scalar feature, and selects the appropriate threshold crossing \cite{ZhaoMaffaSandstede2025}. For this reason, the present algorithm is best described as a \emph{feature-based continuation method inspired by pseudo-arclength continuation}, rather than a full Newton-based pseudo-arclength scheme \cite{ZhaoMaffaSandstede2025}.

This distinction is important for interpretation. The method inherits the flexibility of path-following in parameter space \cite{BeynChampneysDoedelGovaertsKuznetsovSandstede2002,Uecker2021}, but its defining equation is empirical rather than analytical, and its corrector is based on local threshold detection rather than a direct nonlinear solve \cite{ZhaoMaffaSandstede2025}. Accordingly, the continued portions reported in this work should be understood as numerical level sets of simulation-based observables. The lower sweep-based estimates are reported separately and are not assigned the same mathematical status \cite{ZhaoMaffaSandstede2025}.

\section{Spiral Feature Function}

\subsection{Motivation}

For the spiral continuation, the goal is not to continue an exact PDE equilibrium or an explicitly formulated coherent structure. Instead, the aim is to distinguish two qualitatively different late-time regimes of the one-dimensional Brusselator:
\begin{itemize}
    \item a spiral-like or source-defect-like regime, whose spacetime organization is visibly asymmetric, and
    \item a more regular wave/stripe-like regime, whose late-time spacetime behavior is comparatively symmetric.
\end{itemize}
The continuation problem is therefore formulated in terms of a scalar feature that separates these two classes of behavior.

This viewpoint is consistent with the general philosophy of feature-based continuation developed by Zhao, Maffa, and Sandstede, where transition curves are computed from scalar observables extracted from direct simulations rather than solely from a classical boundary-value formulation for exact coherent structures \cite{ZhaoMaffaSandstede2025}. In the present setting, the desired transition curve is represented as a threshold level set of a scalar map
\[
r(\sigma,b),
\]
so that the spiral continuation seeks an approximation of the set
\[
\{(\sigma,b): r(\sigma,b)=r_{\mathrm{thr}}\}.
\]
The central question is therefore not whether an exact analytical bifurcation functional is available, but how to construct a scalar observable that reliably distinguishes spiral/source-defect-like spacetime behavior from more symmetric wave/stripe-like behavior \cite{PerraudBorckmans1993,SandstedeScheel2004,ZhaoMaffaSandstede2025}.

\subsection{PDE Output and Late-Time Observation Window}

The spiral feature is constructed from direct numerical simulations of the PDE. For each parameter pair $(\sigma,b)$ and chosen initial condition $u_0(x)$, the computation is naturally divided into two stages.

First, one performs a \emph{relaxation stage}. Starting from the prescribed initial condition, the PDE is integrated up to a relaxation time $T_{\mathrm{relax}}$. This produces a relaxed state, denoted schematically by
\[
U_{\mathrm{rel}}(x).
\]
In the numerical code used in this work, the full state has two components,
\[
(U_1(x,t),U_2(x,t)),
\]
and the second component is the one used for feature extraction.

Second, starting from the relaxed state, one performs an \emph{observation stage}. The PDE is integrated for an additional observation time interval of length $T_{\mathrm{obs}}$, and the second component is recorded as a spacetime field
\[
V(t,x)=U_2(x,t), \qquad t\in[0,T_{\mathrm{obs}}].
\]
Numerically, this produces a matrix
\[
V \in \mathbb{R}^{N_t \times N_x},
\]
whose rows correspond to time samples and whose columns correspond to spatial grid points. This spacetime array is the raw object from which the spiral detector is constructed.

Because the purpose of the feature is to measure the \emph{late-time} regime rather than transient behavior, the full observation window is not used directly. Instead, one restricts to the final fraction of the recorded spacetime. Writing the total observation field as $V$, define
\[
V_{\mathrm{late}} = V(t,x)\big|_{t\in[t_0,T_{\mathrm{obs}}]},
\]
where
\[
t_0 = (1-\alpha)T_{\mathrm{obs}}
\]
for some chosen fraction $\alpha \in (0,1)$. In the implementation used here, $\alpha=0.25$, so that the feature is computed from the last $25\%$ of the observation window. The spiral detector is therefore a function not of the full transient evolution, but of the restricted late-time spacetime pattern.

\subsection{Symmetry-Defect Construction}

The spiral feature is based on a simple geometric idea: spiral/source-defect-like states are more asymmetric in spacetime than the corresponding wave/stripe-like states \cite{PerraudBorckmans1993,SandstedeScheel2004}. To quantify this, one compares the late-time spacetime field with a transformed copy of itself and measures the mismatch.

In general, let $\mathcal{R}$ denote a reflection operator acting on the spacetime matrix $V_{\mathrm{late}}$. The associated symmetry-defect feature is then defined by
\[
r
=
\frac{\|V_{\mathrm{late}}-\mathcal{R}(V_{\mathrm{late}})\|_F}
{\|V_{\mathrm{late}}\|_F},
\]
where $\|\cdot\|_F$ denotes the Frobenius norm. This is a normalized mismatch:
\begin{itemize}
    \item if the observed late-time spacetime pattern is nearly symmetric under the chosen transformation $\mathcal{R}$, then $r$ is small;
    \item if the pattern is strongly asymmetric with respect to that transformation, then $r$ is larger.
\end{itemize}
Thus, the scalar $r$ serves as a data-driven order parameter for the transition between asymmetric spiral/source-defect-like behavior and more symmetric wave/stripe-like behavior \cite{ZhaoMaffaSandstede2025}.

\begin{figure}[H]
    \centering
    \includegraphics[width=0.8\textwidth]{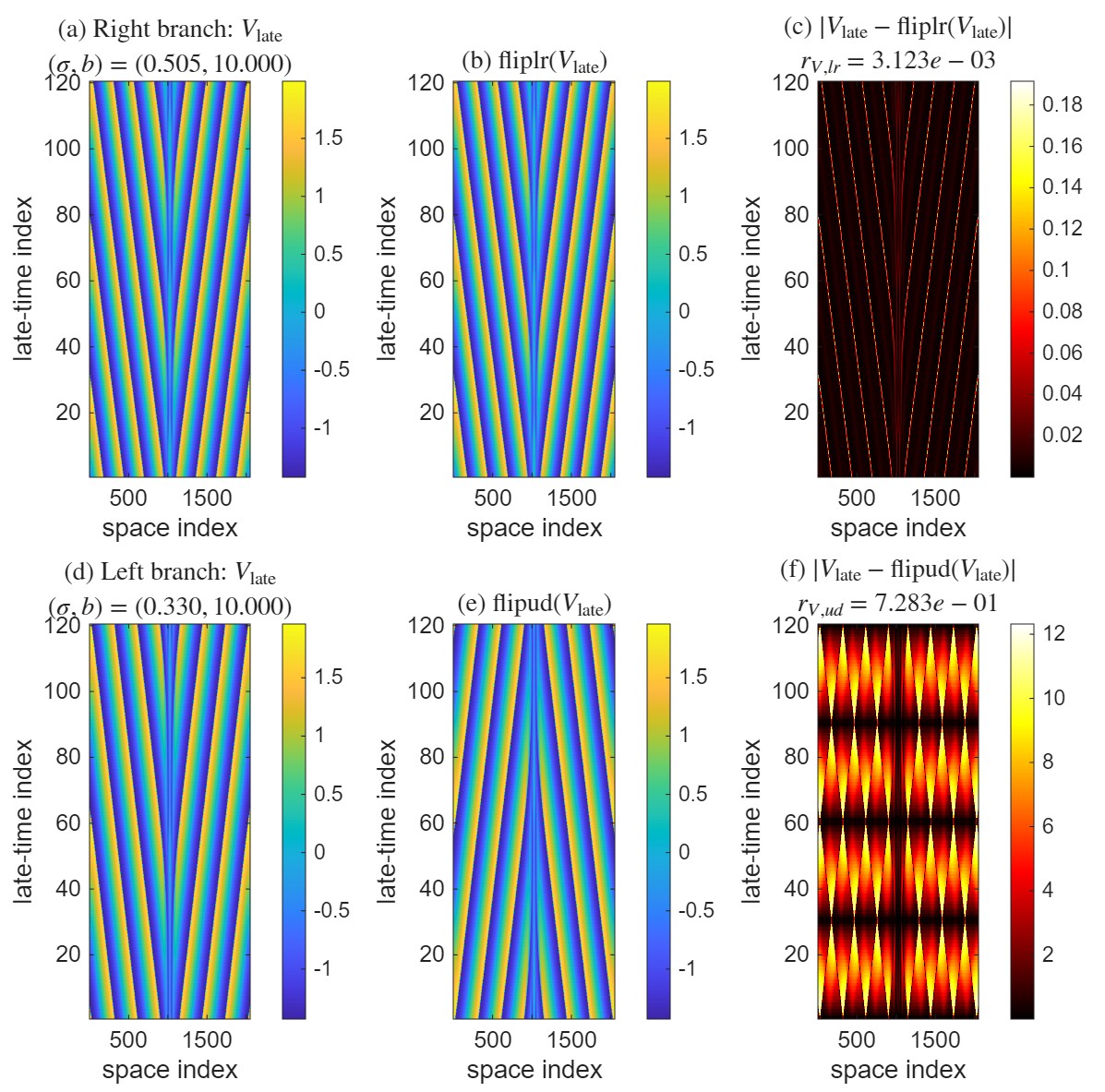}
    \caption{Construction of the spiral symmetry-defect feature from late-time spacetime data. The top row shows a representative right-branch spacetime field $V_{\mathrm{late}}$, its space-reflected copy $\mathrm{fliplr}(V_{\mathrm{late}})$, and the corresponding defect field. The bottom row shows an analogous comparison for a left-branch spacetime field using time reflection $\mathrm{flipud}(V_{\mathrm{late}})$. These comparisons motivate the normalized mismatch metrics used in the spiral feature.}
    \label{fig:spiral_defect_construction}
\end{figure}

The symmetry-defect construction is visualized in Figure~\ref{fig:spiral_defect_construction}, which compares representative late-time spacetime fields with their reflected counterparts. The spiral feature is not based on a subjective visual label, but on a normalized mismatch between the observed late-time spacetime pattern and an appropriately reflected copy of itself. Small mismatch indicates approximate symmetry under the chosen reflection, whereas larger mismatch signals a stronger symmetry defect and hence more pronounced spiral/source-defect-like behavior.

This construction is stronger than a purely snapshot-based detector because it uses the full late-time spacetime field rather than a single terminal profile \cite{PerraudBorckmans1993,SandstedeScheel2004}. In particular, it is sensitive not only to spatial organization but also to temporal phase relations and persistent source-like spacetime structure. As a result, it is better aligned with the actual visual distinction between the numerical regimes studied here.

\subsection{Branch-Dependent Reflection Choice}

Since $V_{\mathrm{late}}$ is a time--space matrix, there are two natural reflection operations available.

The first is \emph{space reflection}, defined by left--right reversal in space,
\[
(\mathcal{R}_x V)(t,x)=V(t,L-x),
\]
which corresponds in matrix language to flipping the columns:
\[
\mathcal{R}_x(V)=\mathrm{fliplr}(V).
\]
The corresponding symmetry-defect feature is
\[
r_{V,\mathrm{lr}}(\sigma,b)
=
\frac{\|V_{\mathrm{late}}-\mathrm{fliplr}(V_{\mathrm{late}})\|_F}
{\|V_{\mathrm{late}}\|_F}.
\]

The second is \emph{time reflection}, defined by reversal in the time direction,
\[
(\mathcal{R}_t V)(t,x)=V(T-t,x),
\]
which in matrix language corresponds to flipping the rows:
\[
\mathcal{R}_t(V)=\mathrm{flipud}(V).
\]
The corresponding feature is
\[
r_{V,\mathrm{ud}}(\sigma,b)
=
\frac{\|V_{\mathrm{late}}-\mathrm{flipud}(V_{\mathrm{late}})\|_F}
{\|V_{\mathrm{late}}\|_F}.
\]

The numerical diagnostics showed that the most reliable detector is not the same on all branches. On right branches, the space-flip detector $r_{V,\mathrm{lr}}$ provided a clean separation between spiral/source-defect-like and wave/stripe-like states, with threshold crossings that aligned well with the actual spacetime transition \cite{PerraudBorckmans1993,ZhaoMaffaSandstede2025}. On left branches, however, the space-flip detector developed large spikes even in stripe-like regions and therefore became unreliable. In that regime, the time-flip detector $r_{V,\mathrm{ud}}$ behaved more robustly and provided a clearer separation of the two pattern classes. This branch dependence is visible in Figure~\ref{fig:spiral_defect_construction}: the reflected comparisons and defect fields illustrate why the space-flip metric is more natural on the right branches, whereas the time-flip metric provides a more robust detector on the left branches.

For this reason, the spiral feature used in continuation is branch-dependent:
\[
r_{\mathrm{spiral}}(\sigma,b)=
\begin{cases}
r_{V,\mathrm{lr}}(\sigma,b), & \text{on right branches},\\[0.5em]
r_{V,\mathrm{ud}}(\sigma,b), & \text{on left branches}.
\end{cases}
\]
It is important to stress that this is an \emph{empirically validated detector choice}, not a universal law of the PDE \cite{ZhaoMaffaSandstede2025}. The branch dependence emerged from numerical diagnostics and was adopted because it gave the most faithful separation of regimes in the specific continuation problem studied here.

\subsection{Comparison with the Earlier Final-Profile Feature}

Before adopting the spacetime-based detector, an earlier and simpler feature based on the terminal spatial profile was also considered. In that approach, if $u_{\mathrm{final}}$ denotes the final profile of the second component, one defines
\[
r_U(\sigma,b)
=
\frac{\|u_{\mathrm{final}}-\mathrm{flip}(u_{\mathrm{final}})\|_2}
{\|u_{\mathrm{final}}\|_2}.
\]
This feature is simpler to compute, since it uses only a single final-time snapshot, but it turned out to be too weak for the present problem.

The difficulty is that two states can have similar terminal profile symmetry while differing substantially in their spacetime organization \cite{PerraudBorckmans1993,SandstedeScheel2004}. In particular, a spiral/source-defect-like state and a more regular wave/stripe state may have comparable final profiles even though their late-time dynamics are qualitatively different. The final-profile detector therefore discards exactly the temporal structure that is most informative for distinguishing the regimes of interest.

The spacetime feature corrects this weakness by replacing the terminal profile with the full late-time spacetime field:
\[
r_U \quad \longrightarrow \quad r_V.
\]
This change is mathematically significant. Instead of measuring symmetry of a single snapshot, one measures symmetry of the entire late-time pattern, including its temporal organization \cite{PerraudBorckmans1993,ZhaoMaffaSandstede2025}. In practice, this yielded a more robust and interpretable detector for the spiral transition.

\subsection{Threshold Definition for Spiral Continuation}

Once the scalar detector has been specified, the spiral transition curve is defined as a threshold level set. Writing the branch-aware feature as $r_{\mathrm{spiral}}(\sigma,b)$, the continued set is
\[
\Gamma_{\mathrm{spiral}}
=
\{(\sigma,b): r_{\mathrm{spiral}}(\sigma,b)=r_{\mathrm{thr}}\}.
\]
Thus, the continuation problem for the spiral regime is to compute an approximation of this feature-defined level set in the $(\sigma,b)$-plane \cite{ZhaoMaffaSandstede2025}.

In this formulation, the feature $r_{\mathrm{spiral}}$ should not be interpreted as an exact analytical bifurcation functional. Rather, it is a \emph{data-driven order parameter} that has been designed to distinguish spiral/source-defect-like asymmetric spacetime behavior from more symmetric wave/stripe behavior \cite{ZhaoMaffaSandstede2025}. Small values of $r_{\mathrm{spiral}}$ indicate that the late-time spacetime pattern is close to the chosen reflected version of itself, while larger values indicate stronger symmetry defect. The continuation uses the empirical observation that across the relevant transition this normalized defect changes sharply enough to support a useful threshold.

Finally, it is worth noting that, because the feature is simulation-based, its value depends not only on $(\sigma,b)$ but also on the seed used to generate the corresponding late-time state \cite{ZhaoMaffaSandstede2025}. This dependence was already important in the continuation framework, where warm-start is helpful along the branch but undesirable within the local detection sweep. The same distinction is relevant here: the intended mathematical object is the feature-defined level set of a stable branch-aware scalar response map, and the use of fixed detection seeds within local sweeps helps make that interpretation as consistent as possible.

\section{Target Feature Function}

\subsection{Motivation from Half-Target / Source-Defect Structure}

The target feature is designed to detect a different class of pattern than the spiral feature. While the spiral feature is based on a spacetime symmetry defect, the target feature is constructed from a more explicitly structural viewpoint. Numerically, the target-like states of interest in this work are associated with half-target or source-defect configurations that exhibit two simultaneous ingredients:
\begin{itemize}
    \item a \emph{structured core} in the interior of the spatial domain, and
    \item an \emph{oscillatory tail} in the far-field region.
\end{itemize}
A state should therefore be classified as target-like only if both of these ingredients are present at the same time \cite{Stich2003,KopellHoward1981}.

This motivates a feature design that differs fundamentally from the spiral case. Instead of measuring a symmetry defect, the target detector measures two separate quantities: one that captures spatial structure in the core, and one that captures temporal oscillation in the tail. The feature is then constructed so that it becomes large only when both of these characteristics are simultaneously present. In this sense, the target detector is a \emph{conjunctive structural detector}: it is intended to recognize the coexistence of a Turing-like core and a Hopf-like tail in the observed late-time spacetime pattern \cite{PerraudBorckmans1993,Stich2003,ZhaoMaffaSandstede2025}.

\subsection{Core and Tail Windows}

As in the spiral case, the target feature is computed from the recorded spacetime field of the second component,
\[
V(t,x)=U_2(x,t),
\]
obtained from a late-time observation window of the PDE simulation. However, rather than using the entire observation field uniformly, the target detector restricts attention both in time and in space.

\begin{figure}[H]
    \centering
    \includegraphics[width=0.85\textwidth]{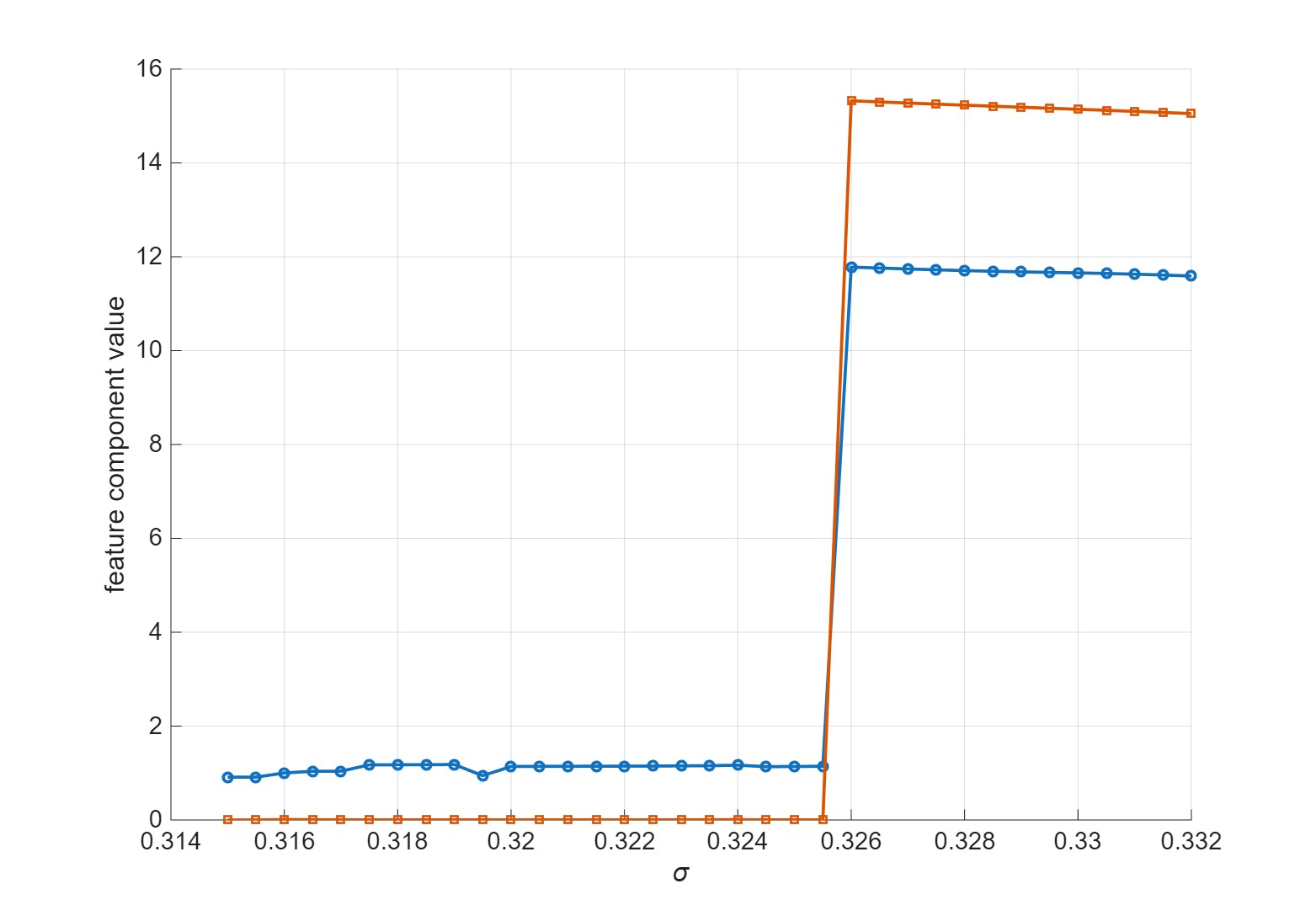}
    \caption{Illustration of the spatial windows used in the target detector for a representative late-time spacetime field. The solid white rectangle marks the core window $I_{\mathrm{core}}$, while the dashed yellow rectangle marks the tail window $I_{\mathrm{tail}}$. These subregions are used to define the core spatial-variance score $S_c$ and the tail temporal-variance score $T_r$.}
    \label{fig:target_windows}
\end{figure}

The spatial subregions used in the target detector are illustrated in Figure~\ref{fig:target_windows}. Figure~\ref{fig:target_windows} shows the geometric decomposition underlying the target feature. The core window is chosen to detect whether the interior of the pattern contains significant spatial structure, while the tail window is chosen to detect whether the far field exhibits sustained temporal oscillation. This spatial separation reflects the intended half-target/source-defect structure of the states being detected.

First, only the final portion of the observation interval is used. If $V \in \mathbb{R}^{N_t \times N_x}$ denotes the recorded spacetime matrix, then the implementation restricts to the last $40\%$ of the observation window:
\[
V_{\mathrm{use}} = V(t,x), \qquad t \geq 0.6T.
\]
Thus, as in the spiral case, the feature is intended to describe the late-time regime rather than transient behavior.

Second, the spatial domain is divided into two subregions:
\begin{itemize}
    \item a \emph{core window}
    \[
    I_{\mathrm{core}}=[0.4L,\,0.6L],
    \]
    centered in the interior of the domain, and
    \item a \emph{tail window}
    \[
    I_{\mathrm{tail}}=[0.8L,\,L],
    \]
    located near the far-right end of the domain.
\end{itemize}
The core window is chosen to measure whether the interior contains significant spatial structure, while the tail window is chosen to measure whether the far field exhibits persistent temporal oscillation. These two subregions reflect the intended qualitative geometry of a half-target/source-defect state \cite{Stich2003,PerraudBorckmans1993}.

\subsection{Core Spatial-Variance Score}

The first ingredient of the target feature is a measure of spatial structure in the core. Let $V_n(x)$ denote the spatial profile at the $n$th retained time sample in the late-time observation window. For each such time slice, compute the spatial variance over the core region $I_{\mathrm{core}}$. Averaging over the retained times gives the score
\[
S_c(\sigma,b)
=
\frac{1}{N_t}
\sum_{n=1}^{N_t}
\operatorname{Var}_{x}\!\bigl(V_n(x)\bigr)\Big|_{x\in I_{\mathrm{core}}}.
\]
Equivalently, in words, $S_c$ is the mean over time of the spatial variance restricted to the core window. As indicated in Figure~\ref{fig:target_windows}, this variance is computed only over the interior core window, so $S_c$ is designed to measure whether the center of the pattern contains persistent spatial structure rather than averaging over the full domain.

The interpretation of this quantity is straightforward. If the interior region is nearly flat or spatially unstructured, then the variance in the core is small and so is $S_c$. If, on the other hand, the interior contains a persistent structured core, then the spatial variance is larger. Thus, $S_c$ acts as a detector for the presence of nontrivial interior patterning \cite{Stich2003,ZhaoMaffaSandstede2025}.

\subsection{Tail Temporal-Variance Score}

The second ingredient is a measure of temporal oscillation in the far-field tail. For each spatial grid point $x_j$ lying in the tail window $I_{\mathrm{tail}}$, compute the temporal variance of the recorded signal $V(t,x_j)$ over the retained late-time interval. Averaging over space gives the score
\[
T_r(\sigma,b)
=
\frac{1}{N_x}
\sum_{j=1}^{N_x}
\operatorname{Var}_{t}\!\bigl(V(t,x_j)\bigr)\Big|_{x_j\in I_{\mathrm{tail}}}.
\]
Equivalently, $T_r$ is the mean over space of the temporal variance restricted to the tail window.

This quantity is intended to detect whether the far field exhibits sustained oscillatory behavior. If the tail is temporally nearly stationary, then the temporal variance is small and so is $T_r$. If the tail contains persistent oscillatory activity, then $T_r$ is larger. Thus, $T_r$ serves as an indicator of Hopf-like temporal behavior in the far field \cite{PerraudBorckmans1993,Stich2003}.

\subsection{Composite Target Feature}

With the core and tail scores in hand, the target feature is defined by
\[
F_{\mathrm{target}}(\sigma,b)
=
\min\{S_c(\sigma,b),\,T_r(\sigma,b)\}.
\]
This is the scalar quantity used in the target continuation.

The use of the minimum is mathematically natural. A half-target/source-defect state should be regarded as target-like only if it has both a structured core and an oscillatory tail \cite{Stich2003,PerraudBorckmans1993}. If one were to combine the two scores by addition, then one very large component could mask the failure of the other. By contrast, taking the minimum enforces a logical ``and'': the feature can be large only when both $S_c$ and $T_r$ are large. If either ingredient disappears, then the minimum drops accordingly. This makes the detector a conjunctive structural observable rather than a generic amplitude score \cite{ZhaoMaffaSandstede2025}.

\begin{figure}[H]
    \centering
    \includegraphics[width=0.7\textwidth]{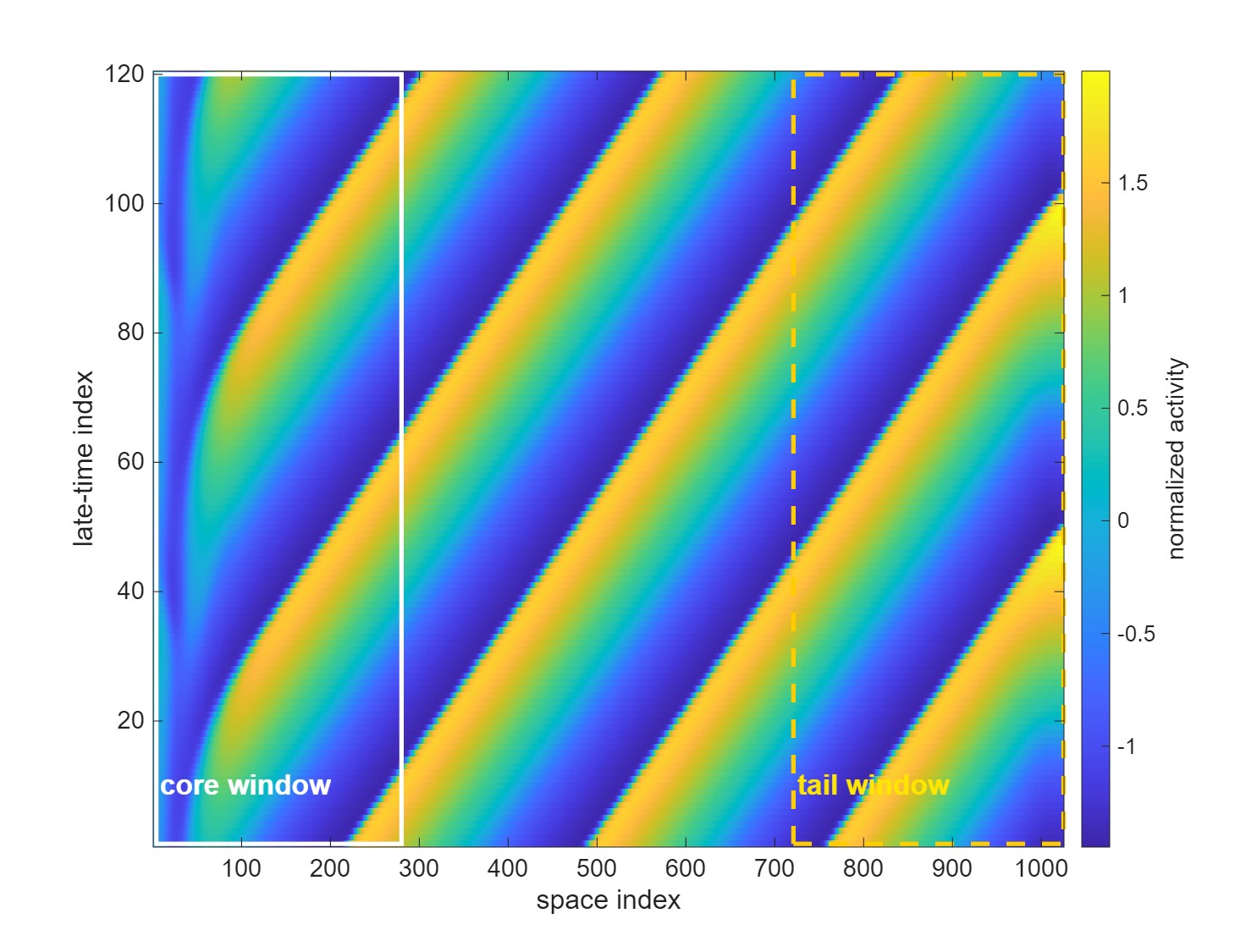}
    \caption{The two components of the target detector along a horizontal sweep at fixed $b=10.02$. The core spatial-variance score $S_c$ and the tail temporal-variance score $T_r$ both increase sharply across the transition region, and the target feature is defined as their minimum. This illustrates why the composite detector acts as a conjunctive structural feature.}
    \label{fig:target_feature_components}
\end{figure}

The behavior of the two target feature components along a representative horizontal sweep is shown in Figure~\ref{fig:target_feature_components}. Figure~\ref{fig:target_feature_components} shows why the minimum is the appropriate way to combine the two component scores. A state should be classified as target-like only when both the core is spatially structured and the tail is temporally oscillatory. Taking the minimum enforces this requirement directly: if either ingredient is absent, the composite feature drops accordingly.

In this sense, the target feature is fundamentally different from the spiral feature of Section~4. The spiral detector is a symmetry-defect observable, whereas the target detector is a conjunctive structural observable \cite{ZhaoMaffaSandstede2025}. The former asks whether the late-time spacetime pattern is asymmetric in an appropriate reflected sense; the latter asks whether two specific dynamical ingredients coexist in the appropriate spatial subregions \cite{PerraudBorckmans1993,Stich2003}.

\subsection{Threshold Definition for Target Continuation}

Once the scalar target detector has been defined, the target transition curve is represented as the threshold level set
\[
\Gamma_{\mathrm{target}}
=
\{(\sigma,b): F_{\mathrm{target}}(\sigma,b)=F_{\mathrm{thr}}\}.
\]
In the computations used here, the preferred threshold is
\[
F_{\mathrm{thr}} = 8.
\]
Accordingly, the continuation framework seeks numerical approximations of the set
\[
\Gamma_{\mathrm{target}}
=
\{(\sigma,b): F_{\mathrm{target}}(\sigma,b)=8\}.
\]

As in the spiral case, this threshold should not be interpreted as an exact analytical bifurcation condition. Rather, $F_{\mathrm{target}}$ is a simulation-derived scalar observable whose level set has been found empirically to separate target-like states from non-target-like states robustly \cite{ZhaoMaffaSandstede2025}. In particular, the threshold value $8$ was chosen because target-like states lie comfortably above it, whereas wave-like or non-target states fall below it. Thus, $F_{\mathrm{target}}$ functions as a data-driven order parameter for the presence of target structure in the one-dimensional Brusselator.

\begin{figure}[H]
    \centering
    \includegraphics[width=\textwidth]{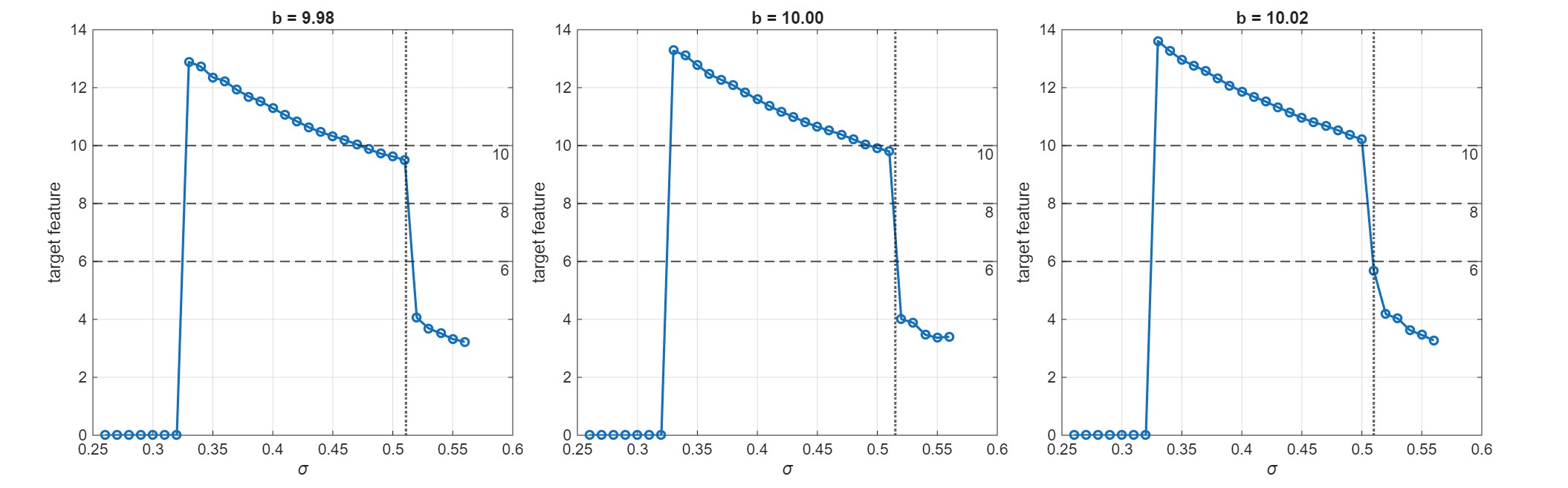}
    \caption{Threshold selection for the target feature from representative horizontal sweeps at fixed $b=9.98$, $b=10.00$, and $b=10.02$. In each case, the target feature exhibits a sharp transition across the sweep, and the level $F_{\mathrm{thr}}=8$ provides a consistent threshold separating target-like states from non-target-like states. The vertical dotted lines indicate the selected transition locations.}
    \label{fig:target_threshold_selection}
\end{figure}

Finally, the target continuation also uses branch-aware crossing logic in practice. In the numerical implementation, the left target branch is tracked by upward threshold crossings of $F_{\mathrm{target}}$ through $8$, while the right target branch is tracked by downward crossings. This is analogous in spirit to the branch-sensitive logic already discussed for the continuation framework: the scalar threshold alone is not enough, and the crossing direction is used to ensure that the correct branch of the target transition set is selected \cite{ZhaoMaffaSandstede2025}. The empirical basis for this threshold choice is illustrated in Figure~\ref{fig:target_threshold_selection}, which shows representative horizontal sweeps of the target feature at nearby values of $b$. Although the feature does not take exactly the same values on every sweep, the level $F_{\mathrm{thr}}=8$ consistently lies in the transition zone and separates the high-feature target regime from the lower-feature non-target regime. This makes $F_{\mathrm{thr}}=8$ a more appropriate continuation threshold than smaller candidate values, since it remains below the clearly target-like states while still lying above the wave-like regime.

With both feature functions now defined, the next section presents the continued side branches, the sweep-based lower transition estimates, and the diagnostics used to distinguish between them.

\section{Numerical Results}

\subsection{Computational Setup}

All numerical transition results reported in this work are based on direct simulation of the one-dimensional Brusselator \eqref{eq:brusselator_1d} together with the feature-based framework developed in the previous sections \cite{Wazwaz2000,ZhaoMaffaSandstede2025}. The calculations are performed in the two-parameter plane $(\sigma,b)$, with all other model parameters fixed as described in Section~2.1. In particular, the reaction parameter $a$ and the second diffusion coefficient $d_2$ are held fixed, while the first diffusion coefficient varies through
\[
d_1=\sigma d_2.
\]

The PDE is solved on a uniform one-dimensional spatial grid with homogeneous Neumann boundary conditions \cite{Pao1982,Murray2003}. Spatial derivatives are approximated by finite differences, producing a semi-discrete system of ordinary differential equations that is integrated in time using \texttt{ode45} in MATLAB \cite{Wazwaz2000,Uecker2021}. The simulation output is recorded as a spacetime array for the second component,
\[
V(t,x)=U_2(x,t),
\]
which serves as the basic data object for the feature functions introduced in Sections~4 and~5.

For each parameter pair $(\sigma,b)$, simulations are divided into a relaxation stage and an observation stage. During relaxation, the system is evolved from a prescribed initial condition so that transient effects decay and the trajectory approaches the relevant late-time regime \cite{ZhaoMaffaSandstede2025}. During observation, the spacetime field of the second component is recorded over a prescribed time interval and then restricted to a late-time window for feature extraction. As discussed earlier, the spiral detector uses a late-time spacetime symmetry-defect construction \cite{PerraudBorckmans1993,SandstedeScheel2004}, while the target detector uses spatial and temporal variance measurements on selected subregions of the spacetime field \cite{Stich2003,ZhaoMaffaSandstede2025}.

The continued side branches are computed by a predictor--corrector procedure \cite{BeynChampneysDoedelGovaertsKuznetsovSandstede2002,ChampneysSandstede2007,ZhaoMaffaSandstede2025}. At each accepted point, a secant predictor advances the computation in the local tangent direction of the previously computed branch \cite{PapakonstantinouTapia2013,AlvesDaSilvaCastroCosta2003}. Correction is then achieved by a local horizontal sweep at fixed $b$, together with threshold detection for the relevant scalar feature. Thus, the method does not solve an explicit augmented nonlinear system, but locates nearby feature-threshold crossings by repeated PDE simulation and feature evaluation \cite{ZhaoMaffaSandstede2025}.

Initialization also requires care. Along an accepted continuation branch, warm-start is used to remain close to the same family of states. Within the local detection sweep, by contrast, fixed seeds or cold-start logic are used whenever possible so that the sampled feature more faithfully represents a well-defined scalar response map rather than a path-dependent artifact of sweep order. In the more delicate lower region, the same fixed seed was used at every value in each vertical $b$-sweep. The transition was then bracketed between the last parameter value showing one clear pattern class and the first neighboring value showing the other. Ambiguous cases were checked through the spacetime plots and were not assigned from the scalar detector alone.

Two scalar thresholds are used for the continued side branches. For the spiral continuation, the relevant threshold is applied to the branch-dependent spacetime symmetry-defect feature defined in Section~4. For the target continuation, the preferred threshold is
\[
F_{\mathrm{thr}}=8,
\]
which empirically separates target-like states from non-target-like states in the side-branch sweeps. These continued portions are interpreted as threshold level sets of simulation-based scalar observables \cite{ZhaoMaffaSandstede2025}. The lower-middle estimates have a different status: they record visually classified transition brackets from vertical sweeps and are not claimed to be points on the same scalar level set.

Finally, the numerical results are not assessed solely from the scalar detector values. At representative points along the continued branches, additional diagnostics are used to verify that the computed curve separates the intended qualitative regimes. These include direct spacetime inspection and fixed-parameter sweeps. In the lower region, direct spacetime classification is part of the boundary-estimation procedure itself, with coherence checks used to identify irregular or mixed patterns \cite{ZhaoMaffaSandstede2025}.

\subsection{Spiral Transition Results}

We begin with the spiral transition results, whose purpose is to detect the boundary between spiral/source-defect-like behavior and the more symmetric wave/stripe-like regime \cite{PerraudBorckmans1993,SandstedeScheel2004,ZhaoMaffaSandstede2025}. As described in Section~4, the continued portions of the spiral transition are defined as the threshold level set of the branch-aware symmetry-defect feature
\[
r_{\mathrm{spiral}}(\sigma,b)=r_{\mathrm{thr}},
\]
where the detector is constructed from the late-time spacetime field of the second component and the reflection operator is chosen according to the branch being continued.

\begin{figure}[H]
    \centering
    \includegraphics[width=0.82\textwidth]{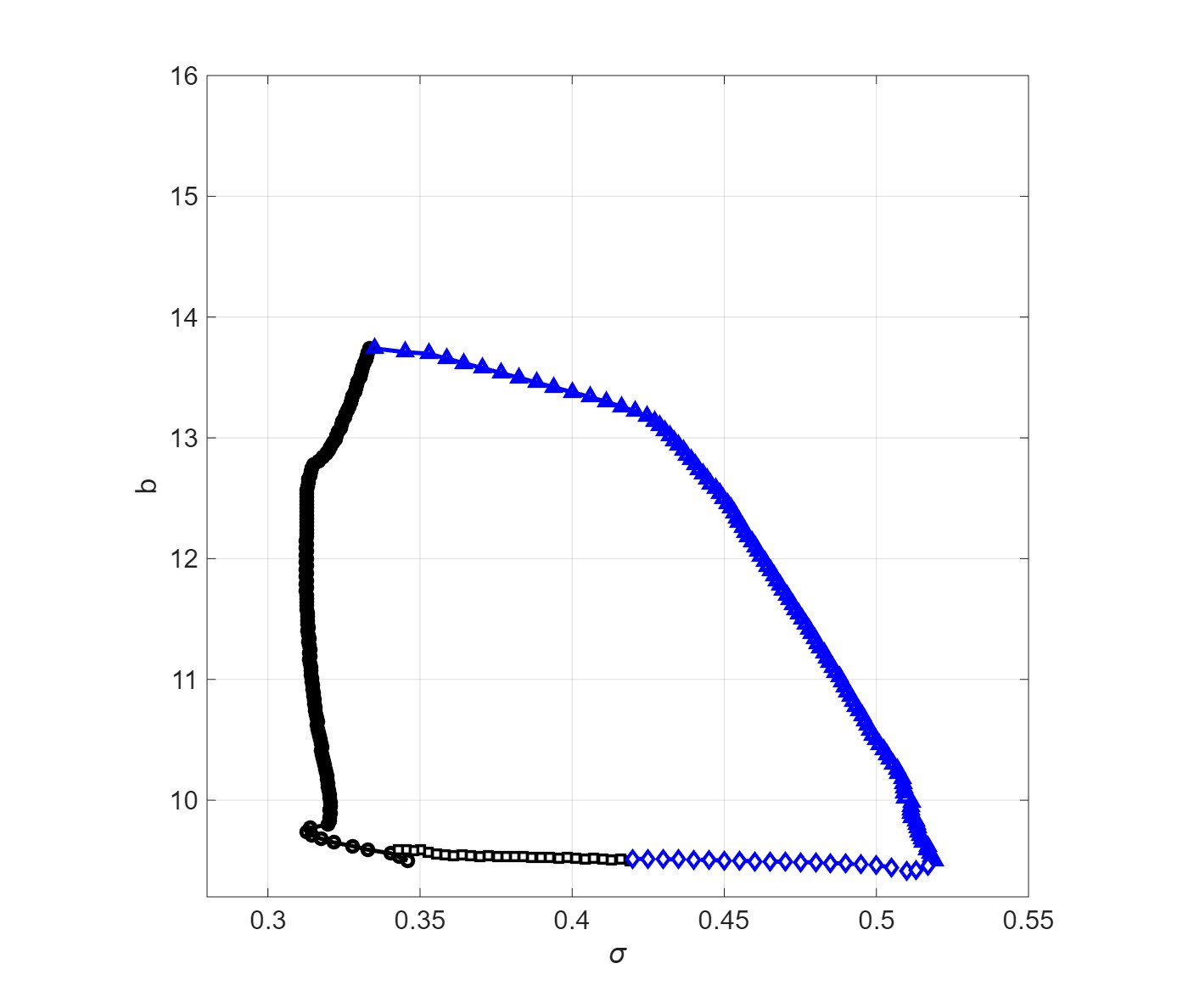}
    \caption{Spiral transition diagram in the $(\sigma,b)$-plane. The upper and side branches are threshold-level sets computed with the horizontal-corrector continuation method. The lower middle-left and middle-right segments are estimates from vertical sweeps in $b$ at selected fixed values of $\sigma$, with the transition identified by direct inspection of late-time spacetime plots. Lines between the sweep-based estimates are included only as visual guides and are not continuation output.}
    \label{fig:spiral_curves_only}
\end{figure}

The resulting spiral transition diagram is shown in Figure~\ref{fig:spiral_curves_only}. Its upper and side portions are recovered by the continuation framework, built around a secant predictor together with a horizontal local corrector \cite{PapakonstantinouTapia2013,AlvesDaSilvaCastroCosta2003,ZhaoMaffaSandstede2025}. In these regions, the geometry of the transition curve is well aligned with
\[
\sigma=\sigma(b),
\]
so fixing $b$ and sweeping in $\sigma$ provides a reliable local correction mechanism.

These side branches form the main left and right walls of the spiral/source-defect region. In the upper and side parameter ranges, the continuation behaves robustly, and the branch-adapted spiral detector gives clean threshold crossings that agree with the visual transition from asymmetric spiral/source-defect-like spacetime to more regular wave/stripe behavior \cite{PerraudBorckmans1993,SandstedeScheel2004}. In this regime, the continuation curves can be interpreted as reliable numerical approximations of the feature-defined regime boundary \cite{ZhaoMaffaSandstede2025}.

A more delicate issue arises in the lower part of parameter space. The original downward continuation suggested that the left-down and right-down branches do not continue independently forever. Instead, the region narrows, and the two lower side branches approach one another, suggesting a lower connecting region. Because the scalar feature becomes less specific there, this part of the boundary was estimated from direct vertical sweeps rather than reported as continuation output.

The lower connecting piece also has a different local geometry from the side branches. The original continuation framework was designed for side-wall geometry, where the curve is naturally viewed as $\sigma=\sigma(b)$ and the corrector fixes $b$ while sweeping in $\sigma$. The lower connection behaves more like
\[
b=b(\sigma),
\]
so the natural local comparison is no longer left/center/right, but up/center/down \cite{BeynChampneysDoedelGovaertsKuznetsovSandstede2002,ChampneysSandstede2007}. At selected values of $\sigma$, $b$ was varied over a discrete grid and the corresponding late-time spacetime plots were compared. The middle-left and middle-right transition estimates were placed between neighboring $b$ values at which the visually observed pattern changed class.

Figure~\ref{fig:spiral_curves_only} therefore combines two types of numerical evidence. The side branches are continued feature thresholds, while the lower middle segments summarize sweep-based visual classifications. Connecting the lower estimates helps display the apparent geometry of the transition region, but the resulting line should not be interpreted as a fully resolved feature-level set.

The lower region remains more delicate than the upper and side branches. A separate vertical-corrector test produced a smooth symmetry-threshold curve that turned downward through the mixed region instead of joining the visually estimated boundary. This behavior indicates that the algorithm continued a detector crossing, but not necessarily the intended transition between coherent spiral/source-defect and wave/stripe patterns. Direct spacetime classification was therefore retained for the reported lower estimates.

Overall, the spiral computations support the following picture. The spiral/source-defect regime has well-resolved upper and side boundaries, while its lower extent is less sharply defined because mixed patterns occur near the apparent transition. The lower sweep-based estimates indicate a connecting geometry, but they should be read as empirical regime-boundary estimates rather than as continued level-set branches \cite{ZhaoMaffaSandstede2025}.

\subsection{Spiral Diagnostics}

The curves and estimates in Figure~\ref{fig:spiral_curves_only} are meaningful only to the extent that they separate the intended qualitative regimes. For this reason, the spiral results were checked by direct spacetime diagnostics \cite{PerraudBorckmans1993,SandstedeScheel2004}. At representative points along the continued side branches, simulations were examined on both sides of the computed curve. In the lower region, the spacetime plots were used directly to assign the transition brackets.

\begin{figure}[H]
    \centering
    \includegraphics[width=0.9\textwidth]{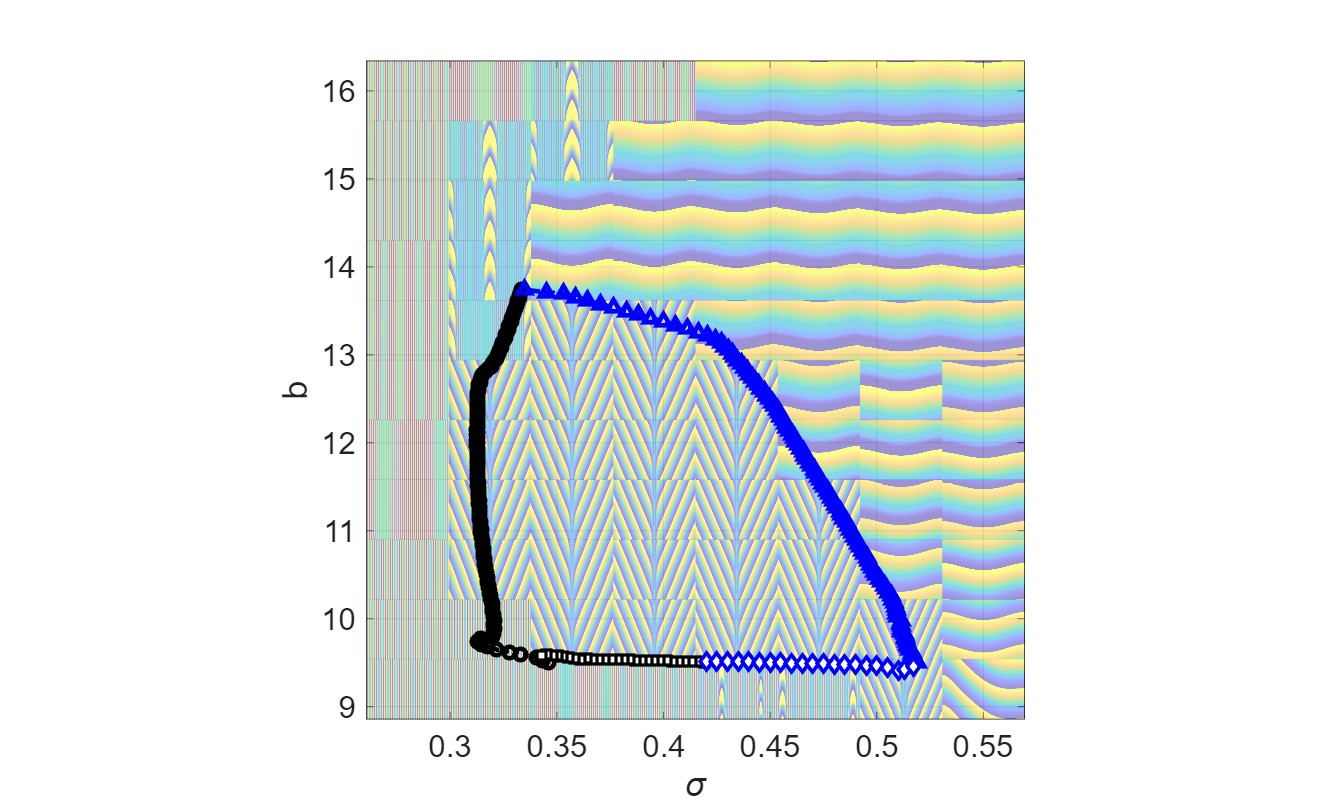}
    \caption{Spiral transition diagram overlaid on a qualitative late-time regime atlas in the $(\sigma,b)$-plane. The background tiles show representative late-time spacetime patterns. The upper and side curves are continued symmetry-threshold branches, while the lower middle segments are sweep-based visual transition estimates. The overlay shows where the scalar continuation is well aligned with the observed regimes and where direct pattern classification is needed.}
    \label{fig:spiral_curves_background}
\end{figure}

Figure~\ref{fig:spiral_curves_background} places the computed spiral branches in their qualitative dynamical context by overlaying them on a background atlas of representative late-time pattern behavior. In the upper and side portions of the boundary, the diagnostics behaved as expected: on one side of the curve, the late-time spacetime field shows the asymmetric source-defect structure that motivates the spiral detector, while on the other side the pattern becomes more regular and wave/stripe-like \cite{PerraudBorckmans1993,CytrynbaumLewis2009}. Thus, in the best-resolved parameter ranges, the continuation is supported not only by threshold crossings of the scalar feature, but also by a direct qualitative change in the observed PDE dynamics \cite{ZhaoMaffaSandstede2025}.

The diagnostics also confirmed that the detector choice had to be branch-dependent. On right branches, the space-reflection defect gave a clean and visually meaningful measure of asymmetry, while on left branches the time-reflection defect proved more robust. This agreement between detector behavior and spacetime inspection shows that the branch-adapted feature choice was not merely a numerical convenience, but a response to genuinely different detector performance across the parameter plane \cite{PerraudBorckmans1993,SandstedeScheel2004}.

The lower connecting region required more careful interpretation. Direct diagnostics suggested that the transition geometry changes qualitatively: instead of continuing as two independent side walls, the lower part of the spiral/source-defect region appears to close through a connecting region. However, mixed or irregular states complicate the interpretation of raw threshold crossings. In such cases, a scalar crossing alone is not sufficient evidence of the intended spiral-to-wave transition, and direct spacetime inspection becomes essential \cite{ZhaoMaffaSandstede2025}.

For this reason, the lower-region diagnostics served two purposes. First, the vertical sweeps located intervals in $b$ over which the visible pattern changed. Second, the spacetime plots and temporal-coherence score helped distinguish clear transitions from detector crossings caused by irregular intermediate states. This distinction also explains why numerical smoothness alone was not used to accept the downward continuation test as the lower spiral boundary.

Taken together, these diagnostics show that the spiral transition diagram cannot be judged solely by the numerical smoothness of a feature-threshold curve. Its validity depends on whether the plotted boundary actually separates the intended pattern classes in the underlying PDE \cite{ChampneysSandstede2007,ZhaoMaffaSandstede2025}. In the upper and side regions, this validation is strong and consistent. In the lower region, the sweep estimates indicate an apparent connecting geometry but also show the increased difficulty of distinguishing coherent regime boundaries from irregular mixed states \cite{ZhaoMaffaSandstede2025}.

\subsection{Target Transition Results}

We next consider the target transition results. As described in Section~5, the target feature is defined by the composite detector
\[
F_{\mathrm{target}}(\sigma,b)=\min\{S_c(\sigma,b),\,T_r(\sigma,b)\},
\]
and the target transition set is represented by the threshold level set
\[
\Gamma_{\mathrm{target}}
=
\{(\sigma,b):F_{\mathrm{target}}(\sigma,b)=F_{\mathrm{thr}}\},
\qquad
F_{\mathrm{thr}}=8.
\]
Thus, continuation of this level set is used to locate the regular side portions of the boundary separating target-like states---characterized by a structured core together with an oscillatory tail---from non-target-like regimes in the $(\sigma,b)$-plane \cite{Stich2003,KopellHoward1981,ZhaoMaffaSandstede2025}.

\begin{figure}[H]
    \centering
    \includegraphics[width=0.75\textwidth]{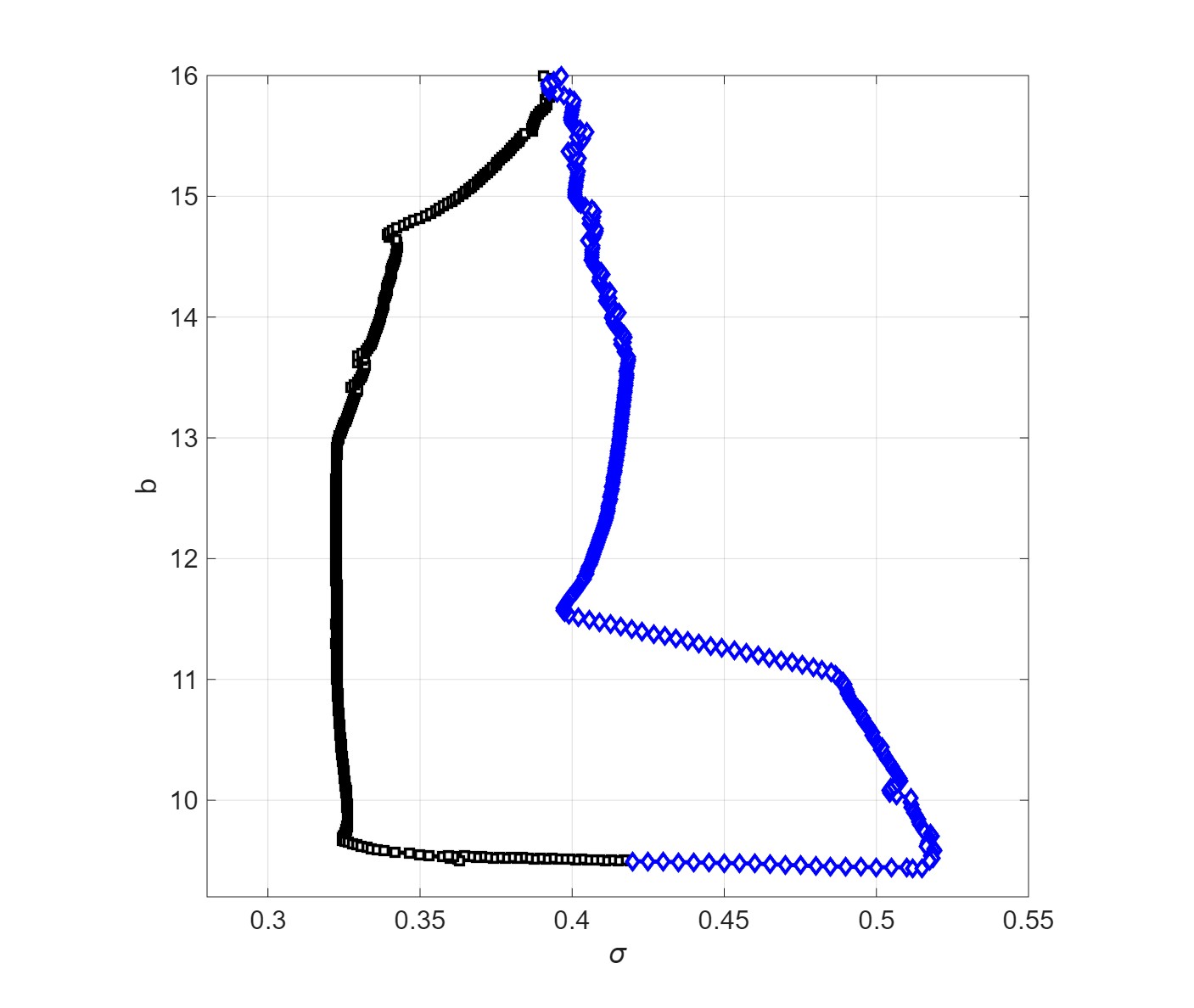}
    \caption{Target transition diagram in the $(\sigma,b)$-plane. The left and right side branches are continued level sets of the target feature. The lower middle-left and middle-right segments are transition estimates obtained from vertical sweeps in $b$ and direct classification of the corresponding late-time spacetime plots. Lines between these sweep-based estimates are visual guides and are not predictor--corrector continuation output.}
    \label{fig:target_curves_only}
\end{figure}

The resulting target transition diagram is shown in Figure~\ref{fig:target_curves_only}. As in the spiral case, the boundary is not represented by a single graph. The figure combines continued left and right side branches with sweep-based middle-left and middle-right estimates in the lower part of parameter space \cite{BeynChampneysDoedelGovaertsKuznetsovSandstede2002,ChampneysSandstede2007,ZhaoMaffaSandstede2025}.

The side branches are the natural output of the original continuation framework. In these regions, the geometry of the target boundary is well aligned with the horizontal corrector, so the transition can be followed efficiently by fixing $b$ and sweeping in $\sigma$ \cite{PapakonstantinouTapia2013,AlvesDaSilvaCastroCosta2003,ZhaoMaffaSandstede2025}. Along these branches, the target feature behaves cleanly: one side of the curve corresponds to states with both a structured core and an oscillatory tail, while the other side corresponds to states for which one or both of these ingredients are lost \cite{Stich2003,ZhaoMaffaSandstede2025}. In this regime, the continuation provides a coherent numerical representation of the target boundary.

As with the spiral transition, the lower part of the parameter plane required a different strategy. The lower side branches suggested that the left and right portions of the target region narrow toward one another, indicating a possible lower connection rather than two side branches extending independently forever. The target feature, however, did not provide a robust threshold crossing for automatic vertical continuation in this region.

The geometric reasoning for the sweeps is the same as in the spiral case. The side branches are naturally treated as portions of a curve of the form
\[
\sigma=\sigma(b),
\]
so that a horizontal corrector is appropriate. By contrast, the lower connecting segment behaves more naturally like
\[
b=b(\sigma),
\]
which calls for the opposite local viewpoint \cite{BeynChampneysDoedelGovaertsKuznetsovSandstede2002,ChampneysSandstede2007,Uecker2021}. At selected fixed values of $\sigma$, $b$ was therefore swept vertically and the late-time spacetime plots were inspected. Transition brackets were assigned where neighboring values of $b$ showed a clear change between the target-like and non-target-like pattern classes. The middle points in Figure~\ref{fig:target_curves_only} summarize these brackets; interpolated portions between sampled estimates serve only to show the apparent geometry.

This distinction is important for interpreting Figure~\ref{fig:target_curves_only}. The lower points support the presence of an approximate cap-like transition region, but they are not evidence that the target level set $F_{\mathrm{target}}=8$ continues smoothly across the entire lower boundary. They instead record where the observed spacetime pattern changed along the sampled vertical slices.

\begin{figure}[H]
    \centering
    \includegraphics[width=0.9\textwidth]{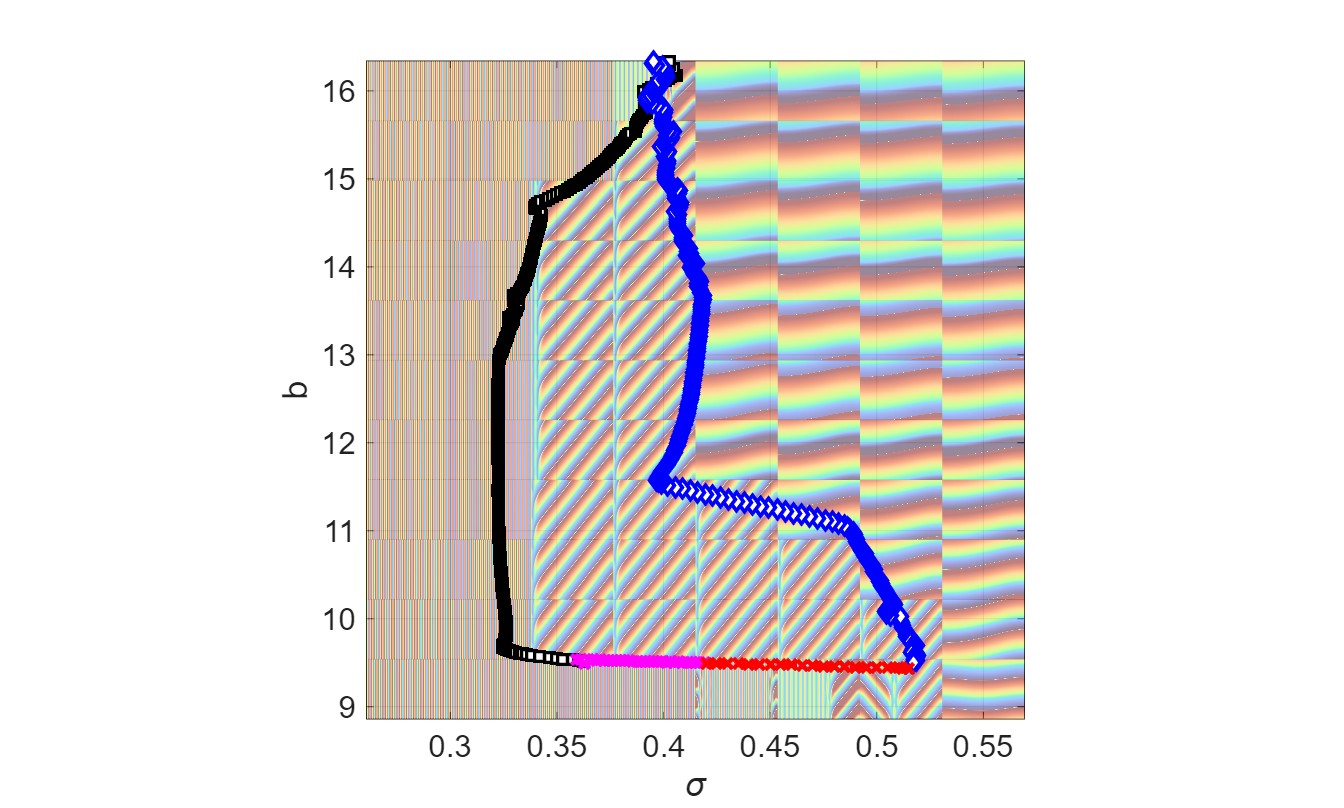}
    \caption{Target transition diagram overlaid on a qualitative late-time regime atlas in the $(\sigma,b)$-plane. The background tiles show representative late-time spacetime patterns at sampled parameter values. The side curves are continued target-feature thresholds, while the lower middle segments are sweep-based visual transition estimates. This visualization places both types of numerical evidence in the broader pattern landscape.}
    \label{fig:target_curves_background}
\end{figure}

Figure~\ref{fig:target_curves_background} places these results in their qualitative dynamical context by overlaying them on a background atlas of representative late-time pattern behavior. The target detector is built from two directly interpretable structural components and gives clean threshold crossings on the side branches \cite{Stich2003,ZhaoMaffaSandstede2025}. The lower region is more difficult. Mixed or irregular states appear near the apparent boundary, and the detector becomes less specific when the system passes through an intermediate regime that is neither cleanly target-like nor cleanly wave-like \cite{ZhaoMaffaSandstede2025}. Direct inspection is therefore part of the lower-boundary estimate rather than only a later validation step.

Overall, the target computations support the view that the target-like regime has well-resolved side boundaries and an approximate lower connecting region. The side branches are resolved by the horizontal-corrector continuation framework, while the lower connection is estimated from vertical sweeps and visual pattern classification. This makes the target problem a useful parallel to the spiral case while also showing the limit of treating every part of a visually observed regime boundary as one continued scalar level set \cite{ZhaoMaffaSandstede2025}.

\subsection{Comparison with Direct Horizontal Sweeps}

An important question for both the spiral and target computations is whether the continued side branches represent genuine regime boundaries or merely internally consistent outputs of the predictor--corrector algorithm \cite{ZhaoMaffaSandstede2025}. For this reason, selected side branches were checked against direct horizontal sweeps at fixed values of $b$. These sweeps provide an independent way to locate transition points by evaluating the relevant scalar feature across a one-parameter slice and identifying the associated threshold crossing \cite{ZhaoMaffaSandstede2025}.

\begin{figure}[H]
    \centering
    \includegraphics[width=0.75\textwidth]{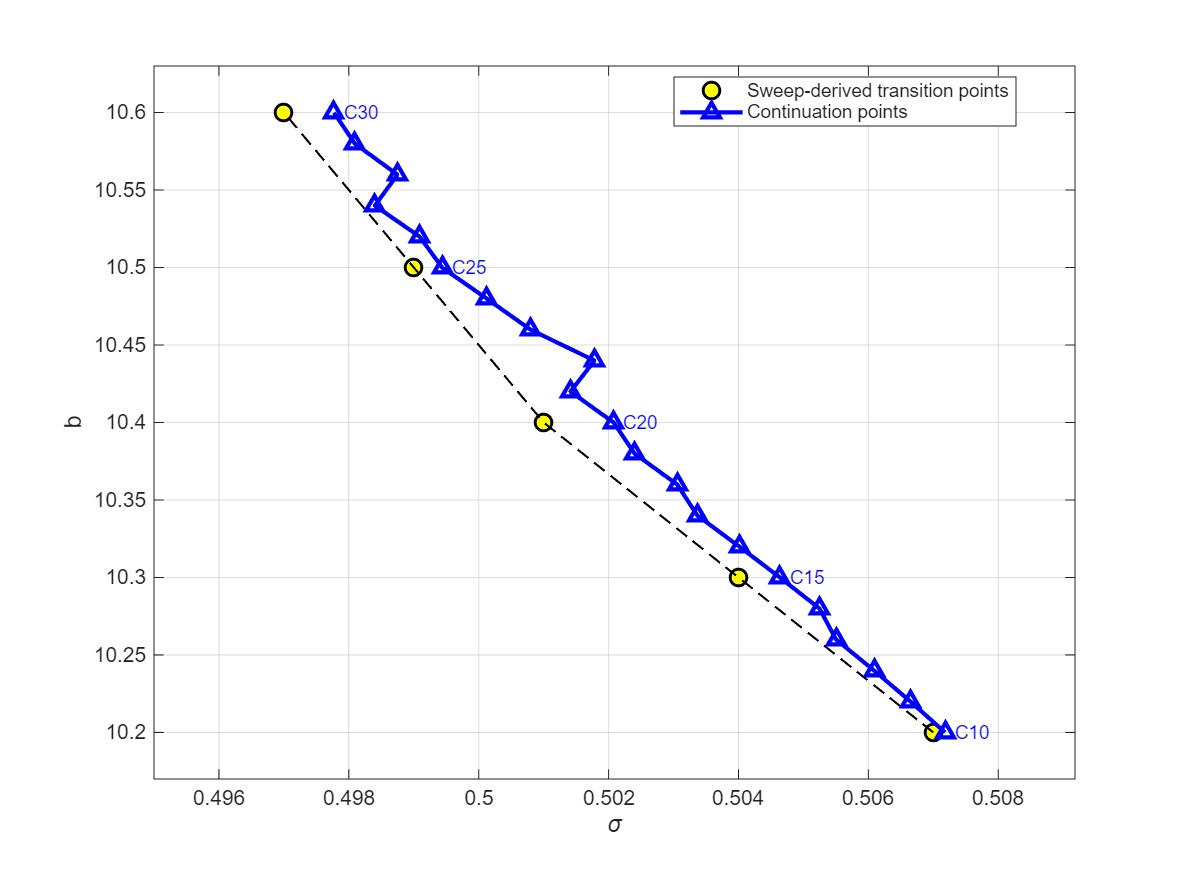}
    \caption{Comparison between a selected continuation segment and independently computed sweep-derived transition points for the target right-up branch. The blue triangles show continuation points along the branch, while the yellow circles show transition points obtained from direct horizontal sweeps at nearby values of $b$. The close alignment of the two sets of points provides an independent validation of the continuation result in this parameter range.}
    \label{fig:sweep_vs_continuation_target}
\end{figure}

The agreement between continuation points and independently detected sweep-derived transition points is illustrated in Figure~\ref{fig:sweep_vs_continuation_target}. In that figure, a selected continuation segment on the target right-up branch is compared with transition points obtained directly from horizontal sweeps over a nearby range of $b$ values. If the continuation is tracking the intended target boundary correctly, then the sweep-derived transition points should lie close to the corresponding continuation segment \cite{ZhaoMaffaSandstede2025}.

Figure~\ref{fig:sweep_vs_continuation_target} shows that the continuation branch is not merely a numerical artifact of the predictor--corrector algorithm. Over the displayed parameter range, the continuation points lie close to the transition locations obtained independently from direct horizontal sweeps. This agreement supports the interpretation of the computed branch as a genuine feature-defined regime boundary rather than as a branch selected only by internal continuation logic \cite{ZhaoMaffaSandstede2025}.

This type of comparison is valuable because it validates the continuation from outside the continuation procedure itself. The continuation algorithm uses local predictors, correctors, hysteresis rules, and branch-selection logic, all of which could in principle bias the computed branch if the detector were poorly behaved \cite{BeynChampneysDoedelGovaertsKuznetsovSandstede2002,ChampneysSandstede2007,ZhaoMaffaSandstede2025}. By contrast, a direct horizontal sweep does not assume the existence of a continuation branch in advance. It simply samples the scalar detector over a one-parameter slice and asks where the threshold crossing occurs. Agreement between these two procedures therefore provides strong evidence that the computed branch is tracking a real feature-defined transition \cite{ZhaoMaffaSandstede2025}.

Such sweep comparisons are inherently local, however. Agreement on a selected side segment does not establish the lower boundary, especially in mixed regions where the detector becomes less specific \cite{ZhaoMaffaSandstede2025}. For this reason, the horizontal comparison should be viewed as one layer of validation for the continued side branches, separate from the vertical sweeps used to estimate the lower region.

In summary, the direct horizontal sweeps provide an important independent check on the numerical continuation results. Where the detector is clean and the continuation geometry is favorable, the sweep-derived transition points closely match the continuation branch. This agreement strengthens the interpretation of the continued side curves as meaningful numerical approximations of regime boundaries in the $(\sigma,b)$-plane \cite{ZhaoMaffaSandstede2025}.

\begin{figure}[t]
    \centering
    \includegraphics[width=\textwidth]{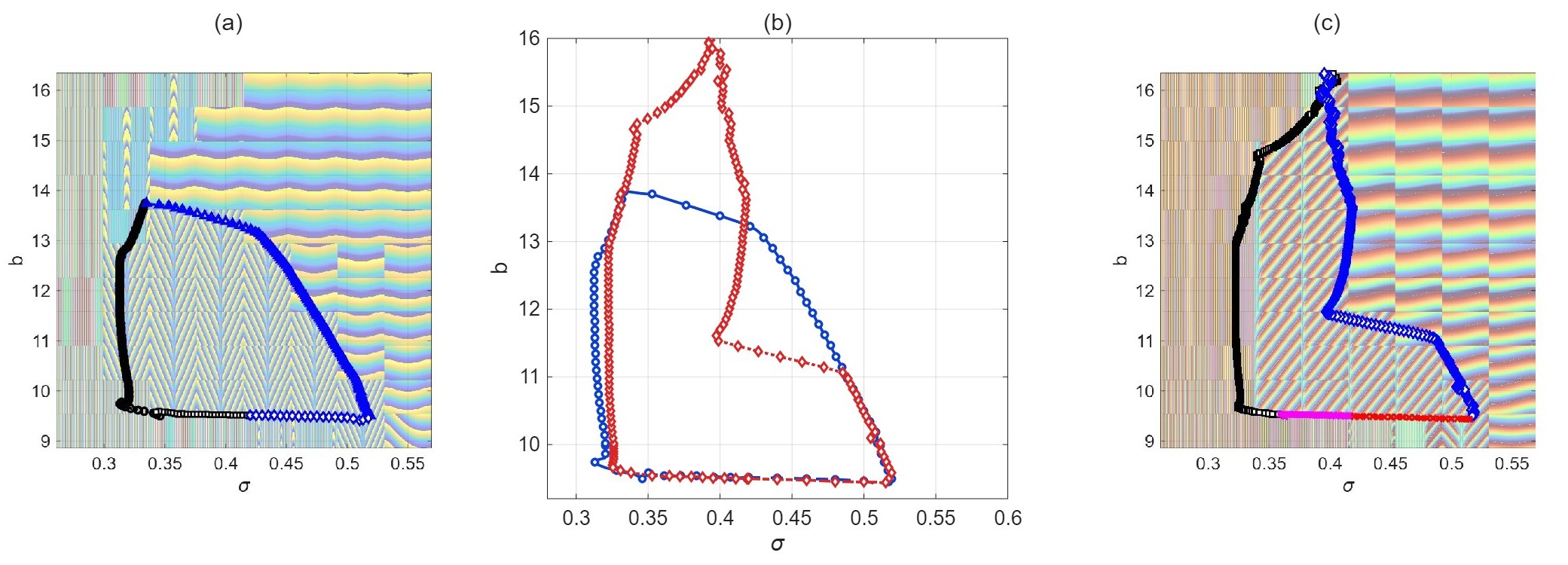}
    \caption{Comparison of the spiral and target transition diagrams in the $(\sigma,b)$-plane. Panel (a) shows the target diagram, panel (b) overlays the spiral and target results using one color for each pattern class, and panel (c) shows the spiral diagram. In each panel, the side branches are continued feature thresholds, whereas the lower middle segments are sweep-based visual transition estimates connected by guides. The panels show both the common large-scale organization and the quantitative differences between the two pattern classes.}
    \label{fig:three_panel_comparison}
\end{figure}

For a direct comparison of the two transition diagrams, Figure~\ref{fig:three_panel_comparison} places the target results, the combined spiral--target overlay, and the spiral results side by side in a common parameter plane.

Figure~\ref{fig:three_panel_comparison} summarizes the overall geometry suggested by the two numerical studies. The spiral and target diagrams have a similar large-scale organization, including well-resolved side branches and approximate lower connecting regions, but their placement and shape are quantitatively distinct. The comparison also separates what is established by feature continuation from what is estimated through direct classification of parameter sweeps.

\subsection{Failure Modes and Ambiguous Regions}

Although the continuation framework performs well over substantial portions of the spiral and target transition sets, there are also parameter regions in which the interpretation of the computed feature becomes less clear \cite{ZhaoMaffaSandstede2025}. These difficulties arise most prominently in lower or extreme parts of parameter space, where the late-time dynamics may no longer fit cleanly into the intended dichotomy between coherent regimes. In such regions, the detector may still cross its nominal threshold, but the crossing need not correspond to the desired transition between well-separated pattern classes \cite{PerraudBorckmans1993,ZhaoMaffaSandstede2025}.

One source of difficulty is geometric. As discussed above, the original continuation framework was designed primarily for side-branch geometry, where the transition curve is locally well represented as $\sigma=\sigma(b)$ and a horizontal corrector is natural. In lower regions, however, the transition set may instead behave more like $b=b(\sigma)$, which motivates a vertical sweep in $b$ \cite{BeynChampneysDoedelGovaertsKuznetsovSandstede2002,ChampneysSandstede2007,Uecker2021}. A vertical corrector can follow scalar threshold crossings in this orientation, but the lower-region tests showed that improved geometric alignment does not ensure that the selected crossing represents the intended physical transition.

A second source of difficulty is dynamical rather than geometric. In some parameter regions, the PDE supports localized-source patterns that are neither cleanly spiral/source-defect-like nor cleanly wave/stripe-like in the sense assumed by the detector design \cite{PerraudBorckmans1993,SandstedeScheel2004,CytrynbaumLewis2009}. Such states may retain a visible localized source or interior core while also exhibiting atypical temporal structure, partial loss of coherence, or mixed spatial organization. As a result, they may still produce scalar detector values suggestive of a transition even though they do not lie on the intended coherent regime boundary \cite{ZhaoMaffaSandstede2025}.

\begin{figure}[H]
    \centering
    \includegraphics[width=\textwidth]{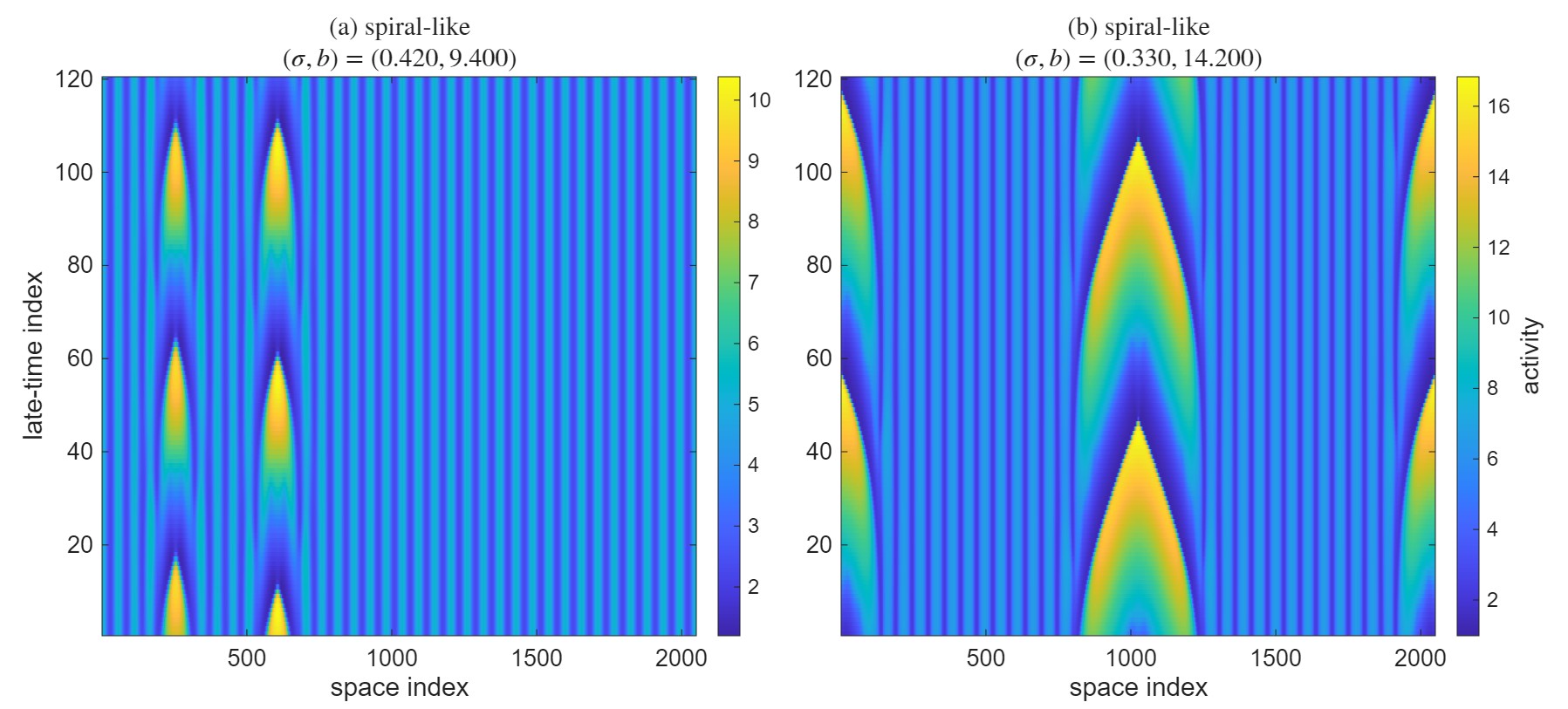}
    \caption{Representative localized-source patterns in ambiguous parameter regions. Although both examples retain source-like localized activity, they do not fit cleanly into the intended regime dichotomy used for continuation. Such states illustrate why threshold crossings alone can become difficult to interpret in delicate parameter regions and why coherence diagnostics and direct spacetime classification are needed.}
    \label{fig:ambiguous_localized_sources}
\end{figure}

Figure~\ref{fig:ambiguous_localized_sources} shows that the difficulty is not merely numerical noise. In both panels, the late-time spacetime field retains a visibly localized source-like structure, yet the resulting pattern is not cleanly aligned with the idealized regime categories used in the continuation framework. States of this kind are important because they can still produce substantial detector values or threshold crossings even though they do not represent the intended coherent transition between well-separated regimes \cite{PerraudBorckmans1993,CytrynbaumLewis2009}. The coherence score can flag some of these states, but the final distinction still requires inspection of the spacetime pattern.

These observations clarify the role of the additional diagnostics introduced in the lower-region computations. The coherence score and direct spacetime inspection respond to a genuine feature-selection problem in the underlying dynamics \cite{ZhaoMaffaSandstede2025}. When the PDE admits patterns that lie in a grey zone between the intended regime categories, no scalar threshold can be interpreted blindly. In that setting, continuation must be supplemented by qualitative validation, and a sweep-based transition bracket may be more honest than an automatically continued feature crossing \cite{ZhaoMaffaSandstede2025}.

For this reason, the continued side branches should be understood as the most reliable portions of a broader regime structure, not as exact analytical boundaries that remain equally sharp in every part of parameter space \cite{ZhaoMaffaSandstede2025}. In well-resolved regions, the detector, continuation geometry, and direct spacetime diagnostics support the same interpretation. In ambiguous regions, detector specificity weakens, mixed states appear, and the lower boundary is reported only as a sweep-based estimate. Recognizing this difference is an important part of assessing the scope of the feature-based continuation framework developed in this work \cite{ChampneysSandstede2007,Uecker2021,ZhaoMaffaSandstede2025}.

\section{Discussion}

\subsection{What Was Achieved}

This work developed and tested a feature-based framework for mapping transitions between qualitatively different late-time regimes in a one-dimensional Brusselator \cite{ZhaoMaffaSandstede2025,Wazwaz2000,YuGumel2001}. The central idea was to replace an explicit analytical bifurcation condition by scalar observables extracted from direct PDE simulations, and then to continue the corresponding threshold level sets where the feature response remained regular \cite{ZhaoMaffaSandstede2025}. In lower mixed regions, the transition was instead estimated through vertical sweeps and direct spacetime classification. The resulting approach combines numerical tracking of feature-defined regime boundaries with a more direct classification method where a single scalar observable is insufficient \cite{ChampneysSandstede2007,Uecker2021,ZhaoMaffaSandstede2025}.

Two distinct feature functions were constructed. The first, used for the spiral transition, was a branch-adapted spacetime symmetry-defect feature computed from the late-time spacetime field of the second component \cite{PerraudBorckmans1993,SandstedeScheel2004}. The second, used for the target transition, was a composite detector defined as the minimum of a core spatial-variance score and a tail temporal-variance score \cite{Stich2003,KopellHoward1981}. These represent two different detector philosophies: the spiral feature measures a normalized symmetry defect in spacetime, whereas the target feature measures the simultaneous presence of two structural ingredients characteristic of a half-target or source-defect state \cite{PerraudBorckmans1993,Stich2003,ZhaoMaffaSandstede2025}.

Using these features together with a secant-based predictor and a local sweep corrector, the work computed nontrivial side branches for both the spiral and target problems \cite{PapakonstantinouTapia2013,AlvesDaSilvaCastroCosta2003,ZhaoMaffaSandstede2025}. The parameter maps also suggested that the regime boundaries are not well described by a single graph over one coordinate \cite{BeynChampneysDoedelGovaertsKuznetsovSandstede2002,ChampneysSandstede2007}. The upper and side branches were recovered reliably by continuation, while lower middle-left and middle-right transition estimates were obtained from vertical sweeps and direct spacetime classification. Together, these results show a two-dimensional regime geometry while preserving a clear distinction between continued feature-level sets and visually estimated portions of the boundary \cite{AllgowerGeorg2000,Uecker2021}.

The numerical study also showed that detector design must be adapted to the pattern class being continued. For the spiral transition, a single universal symmetry detector was not sufficient: the feature had to be branch-dependent, with space reflection performing best on right branches and time reflection on left branches \cite{PerraudBorckmans1993,SandstedeScheel2004}. For the target transition, the continuation succeeded because the detector was tied directly to the intended pattern structure, namely the coexistence of a structured core and an oscillatory tail \cite{Stich2003,KopellHoward1981}. More generally, the success of the continuation depends crucially on whether the chosen scalar feature captures the actual dynamical distinction visible in the PDE \cite{ZhaoMaffaSandstede2025}.

A further contribution of the work was methodological. The computations showed that, in a simulation-based continuation problem, geometric branch-tracking and detector validity must be considered together \cite{ZhaoMaffaSandstede2025}. The lower geometry motivated vertical correction tests because a horizontal corrector was poorly aligned with an apparent boundary of the form $b=b(\sigma)$ \cite{BeynChampneysDoedelGovaertsKuznetsovSandstede2002,ChampneysSandstede2007}. These tests also showed that correct geometry is not enough: ambiguous or irregular patterns can still produce misleading threshold crossings \cite{PerraudBorckmans1993,CytrynbaumLewis2009,ZhaoMaffaSandstede2025}. This motivated branch-aware crossing selection, cold-start local sweeps, coherence diagnostics, and direct pattern classification \cite{ZhaoMaffaSandstede2025}.

Finally, the continuation results were validated by more than the smooth appearance of the computed curves. Direct spacetime diagnostics and fixed-parameter sweeps showed that, in the best-resolved parameter regions, the continued side branches do separate the intended qualitative regimes \cite{PerraudBorckmans1993,ZhaoMaffaSandstede2025}. In the lower region, the same diagnostics exposed the limits of the scalar features and provided the basis for the reported sweep estimates. The features therefore act as meaningful empirical order parameters where their crossings remain specific, but not as universal classifiers across every pattern regime \cite{ZhaoMaffaSandstede2025}.

\subsection{Mathematical Interpretation}

From a mathematical point of view, the continued portions of the transition curves should be interpreted as feature-defined level sets in parameter space \cite{ZhaoMaffaSandstede2025}. More precisely, for a simulation-derived scalar observable $F(\sigma,b)$, the continued object is a set of the form
\[
\Gamma
=
\{(\sigma,b):F(\sigma,b)=F_{\mathrm{thr}}\},
\]
or equivalently,
\[
\Gamma
=
\{(\sigma,b):r(\sigma,b)=r_{\mathrm{thr}}\},
\]
depending on the normalization of the detector. Thus, the continued side branches are not analytical bifurcation curves in the classical sense, but numerical level sets of scalar observables constructed from late-time PDE behavior \cite{ZhaoMaffaSandstede2025}. The lower-middle segments have a different mathematical status: they are empirical transition brackets obtained from one-parameter sweeps and visual classification, with interpolation used only to display their large-scale arrangement.

This distinction is essential. In classical continuation theory, one typically studies solution branches defined as regular zeros of an explicitly formulated nonlinear operator, often augmented by phase conditions and solved by Newton-based boundary-value methods \cite{BeynChampneysDoedelGovaertsKuznetsovSandstede2002,ChampneysSandstede2007,Uecker2021}. By contrast, the present work does not begin with an explicit nonlinear equation whose zero set represents the transition. Instead, the defining scalar is produced by a simulation pipeline,
\[
(\sigma,b)
\longmapsto
\text{PDE solution}
\longmapsto
\text{late-time spacetime data}
\longmapsto
\text{feature value}.
\]
The resulting continuation problem is therefore closer in spirit to the data-driven framework of Zhao, Maffa, and Sandstede, in which pattern transitions are traced through simulation-based observables rather than through a classical boundary-value formulation \cite{ZhaoMaffaSandstede2025}.

This perspective suggests a natural interpretation of the scalar features used here. They function as \emph{data-driven order parameters} for regime transitions \cite{ZhaoMaffaSandstede2025}. The spiral feature measures a branch-adapted symmetry defect in late-time spacetime \cite{PerraudBorckmans1993,SandstedeScheel2004}; the target feature measures the joint presence of a structured core and an oscillatory tail \cite{Stich2003,KopellHoward1981}. In this respect, the spiral detector is also loosely connected to the broader mathematical theme of symmetry breaking: it is not a theorem-level symmetry-breaking functional, but it is designed to quantify the emergence of a less symmetric regime from a more symmetric one in a numerically tractable way \cite{SmollerWasserman1990}. Neither feature is a rigorous analytical normal form coefficient or instability functional, yet each is designed to respond sharply enough across the transition to serve as a usable detector \cite{ZhaoMaffaSandstede2025}.

At the same time, the mathematical meaning of such a feature depends on how it is sampled. Because the detector is produced by direct simulation, it is more accurate to regard it as
\[
F(\sigma,b;u_0),
\]
with an implicit dependence on the initial condition or seed used to generate the late-time state. This is why the distinction between warm-start and cold-start matters mathematically. Along a continuation branch, warm-start is useful because it helps follow a consistent family of states. Within a local sweep, however, warm-start introduces path dependence and weakens the interpretation of the detector as a well-defined scalar response map over the swept parameter \cite{ZhaoMaffaSandstede2025}. The use of fixed detection seeds in local sweeps therefore makes the continuation problem closer to the mathematical ideal of tracking a genuine level set of a function of $(\sigma,b)$.

The work also suggests a geometric interpretation of the different pieces seen numerically. The continued side branches are locally well represented as $\sigma=\sigma(b)$, whereas the sweep-based lower transition estimates are arranged more naturally as $b=b(\sigma)$ \cite{BeynChampneysDoedelGovaertsKuznetsovSandstede2002,ChampneysSandstede2007}. The regime boundary should therefore not be thought of as a single graph. At the same time, the lower estimates do not establish that all of these pieces belong to one regular scalar level set \cite{AllgowerGeorg2000,Uecker2021}.

More broadly, the results support the idea that meaningful dynamical information can be extracted even when exact analytical bifurcation curves are unavailable \cite{ZhaoMaffaSandstede2025}. In the present problem, the continued branches and sweep-based estimates organize the observed pattern landscape in a useful way \cite{CrossHohenberg1993,ZhaoMaffaSandstede2025}. They identify where detector-defined transitions are reliable, suggest how the lower regime boundaries are arranged, and show which parts of parameter space remain ambiguous.

\subsection{Strengths and Limitations}

A major strength of the framework developed in this work is its flexibility. The method does not require that the pattern of interest be available as an explicitly formulated equilibrium, traveling wave, or coherent structure satisfying a known boundary-value problem \cite{ChampneysSandstede2007,Uecker2021}. Instead, it works directly with the output of time-dependent PDE simulations \cite{ZhaoMaffaSandstede2025}. This makes it useful in settings where direct simulation reveals robust qualitative regimes, but where the appropriate analytical formulation for classical continuation is difficult to derive in advance \cite{ZhaoMaffaSandstede2025}.

A second strength is that the framework is feature-modular. The same predictor--corrector architecture was used successfully with two quite different detectors: a branch-adapted spacetime symmetry-defect feature for the spiral transition and a conjunctive structural feature for the target transition \cite{PerraudBorckmans1993,Stich2003,ZhaoMaffaSandstede2025}. This suggests that the continuation method itself is not tied to a single scalar observable, but can support different detectors tailored to different pattern classes.

A third strength is that the continuation was validated independently. The continued curves were not accepted simply because the branch-tracking appeared smooth. Instead, direct spacetime diagnostics and fixed-parameter sweeps were used to verify that, in well-resolved regions, the curves do correspond to genuine transitions in the underlying PDE dynamics \cite{PerraudBorckmans1993,ZhaoMaffaSandstede2025}. In a data-driven continuation setting, this combination of continuation and external validation is especially important \cite{ZhaoMaffaSandstede2025}.

At the same time, the work revealed several clear limitations. The first is detector dependence. The continuation is only as good as the scalar feature being continued, and the computations showed that a poorly chosen detector can fail even when the continuation algorithm itself is working properly \cite{ZhaoMaffaSandstede2025}. This was evident in the spiral problem, where a final-profile symmetry metric was too weak and even the spacetime symmetry detector had to be made branch-dependent \cite{PerraudBorckmans1993,SandstedeScheel2004}. It also appeared in lower ambiguous regions, where threshold crossings could occur in mixed or irregular states that were not part of the intended coherent transition \cite{CytrynbaumLewis2009,ZhaoMaffaSandstede2025}.

A second limitation is threshold dependence. The continued branches are defined by level sets of simulation-derived observables, so their location depends on the selected threshold values \cite{ZhaoMaffaSandstede2025}. Although the thresholds used here were chosen empirically and checked by representative sweeps, they are not uniquely determined by theory. The continued branches should therefore be interpreted as threshold-defined regime boundaries rather than exact universal bifurcation sets. The sweep-based lower estimates reduce dependence on one threshold, but introduce uncertainty from the discrete parameter grid and visual pattern classification \cite{ZhaoMaffaSandstede2025}.

A third limitation is geometric. The continuation framework uses a one-dimensional local corrector, which makes the method simpler than a full two-dimensional Newton corrector but also more sensitive to local curve orientation \cite{BeynChampneysDoedelGovaertsKuznetsovSandstede2002,ChampneysSandstede2007,Uecker2021}. The lower-region tests made this especially clear. Changing from a horizontal to a vertical corrector improved alignment with the apparent curve, but did not resolve the separate problem that the detector could select a crossing inside a mixed regime.

A fourth limitation is that the full transition diagram is not fully autonomous. In difficult regions, additional judgment enters through hysteresis windows, branch-aware crossing rules, coherence diagnostics, and visual classification of spacetime plots \cite{ZhaoMaffaSandstede2025}. These ingredients are justified by the PDE behavior, but they also show that the present method is not yet a push-button algorithm.

Finally, some limitations come from the dynamics itself rather than from the algorithm. In ambiguous regions, the PDE supports patterns that do not fit cleanly into the regime categories used to design the detectors \cite{PerraudBorckmans1993,CytrynbaumLewis2009}. In such cases, the ambiguity is not merely numerical; it reflects a genuine grey zone in the pattern dynamics. No scalar detector can remove that ambiguity completely without stronger modeling assumptions \cite{ZhaoMaffaSandstede2025}.

Overall, the framework is most effective when three conditions hold simultaneously: the detector is specific to the transition of interest, the local branch geometry is compatible with the chosen corrector, and the target regimes are genuinely separated in the PDE dynamics. When these conditions hold, the method can recover meaningful regime boundaries efficiently \cite{ZhaoMaffaSandstede2025}. When they fail, continuation becomes correspondingly more delicate and requires additional diagnostics \cite{ChampneysSandstede2007,Uecker2021,ZhaoMaffaSandstede2025}.

\subsection{Future Directions}

The results of this work suggest several natural directions for future work, both on the numerical methodology and on the underlying pattern dynamics.

A first direction is to improve the continuation framework itself by replacing the present one-dimensional local corrector with a more intrinsic two-dimensional correction scheme \cite{BeynChampneysDoedelGovaertsKuznetsovSandstede2002,ChampneysSandstede2007,Uecker2021}. The horizontal corrector recovered substantial portions of the spiral and target side boundaries, while the vertical-corrector tests remained sensitive both to curve geometry and to detector ambiguity. A natural next step would be an adaptive corrector that chooses its search direction from the local tangent or normal geometry, or more ambitiously, a full two-dimensional Newton-type corrector for feature-defined level sets.

A second direction is to refine the feature functions themselves. The present work showed clearly that continuation quality depends strongly on detector quality \cite{ZhaoMaffaSandstede2025}. For the spiral transition, the branch-adapted spacetime symmetry-defect detector worked much better than the earlier final-profile metric, but it was still not perfectly specific in delicate lower regions \cite{PerraudBorckmans1993,SandstedeScheel2004}. For the target transition, the minimum-based structural detector was effective because it encoded the simultaneous presence of a structured core and an oscillatory tail, but other feature constructions may sharpen the boundary further or reduce sensitivity to ambiguous intermediate states \cite{Stich2003,ZhaoMaffaSandstede2025}.

A third direction is to study the ambiguous and lower-region dynamics more systematically rather than treating them only as obstacles to continuation. In the present work, such states were important mainly because they complicated threshold detection and required direct spacetime classification. However, they may also be interesting dynamical objects in their own right \cite{PerraudBorckmans1993,CytrynbaumLewis2009}. The appearance of localized-source patterns that do not fit cleanly into the idealized spiral/source-defect or wave/stripe classes suggests that the PDE may support a richer hierarchy of intermediate regimes than the current detector language captures.

A fourth direction is to connect the numerical continuation more closely with analytical bifurcation theory. The motivating literature, especially the one-dimensional ``spirals'' of Perraud \emph{et al.}\ and the broader pattern-formation framework of Cross and Hohenberg, suggests that the relevant structures may be related to interactions between oscillatory and spatial instabilities \cite{PerraudBorckmans1993,CrossHohenberg1993}. The present work did not attempt a full analytical reduction of the one-dimensional Brusselator in the parameter range studied, but the computed regime boundaries may provide useful guidance for such an analysis. One could ask, for example, whether the feature-defined transition curves correlate with analytically identifiable instability thresholds, codimension-two interactions, or reduced-amplitude-equation predictions \cite{TzouMaBaylissMatkowskyVolpert2013,YuGumel2001}.

A fifth direction is to extend the framework to other pattern-forming PDEs and other classes of transitions. One of the main conceptual advantages of the present approach is that it is detector-based rather than model-specific \cite{ZhaoMaffaSandstede2025}. In principle, the same predictor--corrector architecture could be used in other reaction--diffusion systems, provided one can design a scalar observable that separates the regimes of interest.

Finally, there is a natural computational direction involving automation of validation. In this work, direct spacetime inspection, sweep comparison, and coherence diagnostics were all essential in determining which branches were trustworthy. A useful next step would be to incorporate more of this validation directly into the continuation logic, for example by tracking detector confidence, regularity measures, or sweep consistency alongside the feature value itself \cite{ZhaoMaffaSandstede2025}. This could produce not only a candidate transition curve, but also a local reliability score.

Taken together, these directions suggest that the present work is best viewed not as the endpoint of a completed continuation theory, but as a first framework for studying simulation-defined pattern transitions in reaction--diffusion systems. The main challenge for future work is to preserve the flexibility of the present approach while improving its geometric robustness, detector specificity, and analytical interpretability \cite{ChampneysSandstede2007,Uecker2021,ZhaoMaffaSandstede2025}.

\section{Conclusion}

This work studied data-driven transition mapping for pattern changes in a one-dimensional Brusselator \cite{Wazwaz2000,YuGumel2001,ZhaoMaffaSandstede2025}. The central problem was to determine how qualitatively different late-time regimes are organized in the two-parameter plane $(\sigma,b)$ and to estimate the curves separating them \cite{CrossHohenberg1993,ZhaoMaffaSandstede2025}. Rather than relying on a classical boundary-value formulation for exact coherent structures, the work extracted scalar observables from direct PDE simulations and continued their threshold level sets in parameter regions where the resulting crossings were regular \cite{ZhaoMaffaSandstede2025}.

Two different feature constructions were developed. For the spiral transition, the relevant scalar quantity was a branch-adapted late-time spacetime symmetry-defect feature designed to distinguish asymmetric spiral/source-defect-like behavior from more symmetric wave/stripe behavior \cite{PerraudBorckmans1993,SandstedeScheel2004}. For the target transition, the feature was a composite detector defined as the minimum of a core spatial-variance score and a tail temporal-variance score, thereby enforcing the simultaneous presence of the two ingredients characteristic of a target-like or half-target/source-defect state \cite{Stich2003,KopellHoward1981}. These two detectors reflected two distinct philosophies of feature design, but both fit naturally into the same continuation architecture \cite{ZhaoMaffaSandstede2025}.

The continuation framework combined a secant-based predictor with a local sweep corrector, together with branch-aware crossing rules, hysteresis, and careful treatment of warm-start versus cold-start logic \cite{PapakonstantinouTapia2013,AlvesDaSilvaCastroCosta2003,ZhaoMaffaSandstede2025}. In the side-branch geometry, this framework successfully recovered upper and side portions of the spiral and target transition sets. In the lower part of parameter space, vertical sweeps and direct inspection of spacetime plots were used to estimate middle-left and middle-right transition locations. Tests with a vertical corrector showed that a smooth scalar threshold crossing could continue into mixed pattern regions and therefore should not automatically be interpreted as the intended regime boundary \cite{BeynChampneysDoedelGovaertsKuznetsovSandstede2002,ChampneysSandstede2007,Uecker2021,ZhaoMaffaSandstede2025}.

The numerical results showed that the continued side branches are meaningful in the best-resolved parameter ranges. Direct spacetime diagnostics and comparison with independent horizontal sweeps confirmed that these branches separate the intended qualitative regimes rather than merely reflecting the internal logic of the continuation algorithm \cite{PerraudBorckmans1993,ZhaoMaffaSandstede2025}. The lower sweep estimates have a more limited interpretation: they record visually identified transition brackets and the lines between them are only guides. Neither type of curve is an exact analytical bifurcation set. The continued portions are better understood as numerical level sets of carefully constructed simulation-based observables, while the lower estimates summarize direct classification of the simulated patterns \cite{ChampneysSandstede2007,ZhaoMaffaSandstede2025}.

In this sense, the main contribution of the work is both methodological and conceptual. Methodologically, it demonstrates that predictor--corrector continuation can be adapted to feature-defined transition problems even when no explicit boundary-value formulation is available \cite{ChampneysSandstede2007,Uecker2021,ZhaoMaffaSandstede2025}. Conceptually, it shows that regime geometry can be studied by combining feature continuation in regular regions with direct sweep-based classification where the dynamics are mixed, provided that the different levels of numerical evidence are stated clearly \cite{CrossHohenberg1993,ZhaoMaffaSandstede2025}.

More broadly, the work suggests that there is a productive middle ground between coarse direct simulation and fully classical coherent-structure continuation \cite{ChampneysSandstede2007,ZhaoMaffaSandstede2025}. Feature-based continuation does not replace analytical bifurcation theory, but it provides a useful framework for organizing and studying pattern transitions when the relevant regimes are visible numerically before they are fully understood analytically \cite{CrossHohenberg1993,ZhaoMaffaSandstede2025}. In that sense, the work offers a practical approach to exploring the geometry of pattern selection in reaction--diffusion systems and a starting point for future refinements that may connect simulation-defined transition curves more closely with the underlying theory.

\section*{Code Availability}

The MATLAB code, required seed files, continuation data, CCV job scripts, and figure-generation scripts used in this work are available in the accompanying GitHub repository \cite{YuBrusselatorRepository2026}.

\section*{Acknowledgments}

The author thanks Bj{\"o}rn Sandstede for his guidance throughout this project and for many helpful discussions.

\bibliographystyle{plainnat}
\bibliography{references}

\end{document}